\documentclass{amsart}
\usepackage{amsmath}
\usepackage{amssymb}
\usepackage[all]{xy}

\usepackage[colorlinks=true,linkcolor=black,citecolor=black]{hyperref}

\theoremstyle{plain}
\newtheorem{thm}{Theorem}[section]
\newtheorem{thm*}{Theorem}[section] 
\newtheorem{cor}[thm]{Corollary}
\newtheorem{defn}[thm]{Definition}
\newtheorem{prop}[thm]{Proposition}
\newtheorem{lemma}[thm]{Lemma}

\newtheorem{claim*}{Claim}

\newtheorem{rem}[thm]{Remark}

\numberwithin{equation}{thm}

\newcommand{\bE}{\mathbf{Ext}}
\newcommand{\bR}{{\mathbf R}}
\newcommand{\bu}{\bullet}
\newcommand{\V}{\mathcal V}
\newcommand{\B}{{\mathcal P^{op} \times \mathcal P}}
\newcommand{\biF}{{\mathcal F^{op} \times \mathcal F}}

\newcommand{\F}{\mathcal F}
\newcommand{\Po}{\mathcal P}
\newcommand{\h}{\mathcal H}
\newcommand{\Z}{\mathbb Z}
\newcommand{\bF}{\mathbb F}

\newcommand{\Hom}{\mathrm{Hom} }

\newcommand{\hc}{\mathcal{H}om }
\newcommand{\Ext}{\mathrm{Ext} }
\newcommand{\coh}{\mathrm{H} }
\newcommand{\GL}{\mathrm{GL} }
\newcommand{\gl}{{g\ell}}
\newcommand{\g}{\hc (I,I)}
\newcommand{\K}{K}
\newcommand{\Om}{\Omega }
\newcommand{\Sym}{\mathfrak{S}}
\newcommand{\rk}{\mathrm{rk} }

\newcommand{\I}{{\mathrm I}}
\newcommand{\II}{\mathrm{II}}

\newcommand{\ol}{\overline}

\begin{document}

\title[Cohomology of Bifunctors\ -- \ \today]{Cohomology of Bifunctors}

\author[Vincent Franjou \& Eric M. Friedlander\ -- \ \today]{Vincent
Franjou$^*$ and Eric M. Friedlander$^{**}$}

\address{Laboratoire Jean Leray, Universit\'e de Nantes, BP 92208 -
44322 Nantes Cedex 3, France}
\email{vincent.franjou@univ-nantes.fr}

\address {Department of Mathematics, Northwestern University,
Evanston, IL  60208} \email{eric@math.northwestern.edu}

\date{\today}

\thanks{$^{*}$ The first author is partially supported by the LMJL -
Laboratoire de Math\'ematiques Jean Leray,
CNRS: Universit\'e de Nantes, \'Ecole Centrale de Nantes}
\thanks{$^{**}$ The second author is partially supported by the NSF}


\begin{abstract}
We initiate the study of the cohomology of (strict polynomial)
bifunctors by introducing the foundational formalism, establishing
numerous properties in analogy with the cohomology of functors, and
providing computational techniques.  Since one of the initial
motivations for the study of functor cohomology was the
determination of $\coh^*(\GL(k),S^*(\gl) \otimes \Lambda^*(\gl))$,
we keep this challenging example in mind as we achieve numerous
computations which illustrate our methods.
% Put abstract here
\end{abstract}

\maketitle

% \tableofcontents

%%%%%%%%%%%%%%%%%%%%%%%%
        %Section 0
%%%%%%%%%%%%%%%%%%%%%%%%%%

\section{Introduction}

We fix a prime $p$, a base field $k$ of characteristic $p$, and
consider the category $\V$ of finite dimensional $k$-vector spaces
and $k$-linear maps. The study of the cohomology of categories of
functors from $\V$ to $k$-vector spaces has had numerous
applications, including insight into the structure of modules for
the Steenrod algebra \cite{S} and proof of finite generation of
the cohomology of finite group schemes \cite{FS}.  The computational
power of functor cohomology arises as follows:  the abelian category
of strict polynomial functors $\Po$ of bounded degree enjoys many
pleasing properties which lead to various cohomological computations
(cf. \cite{FFSS});  this cohomology for the category $\Po$ is
closely related to the cohomology for the abelian category of all
functors $\F$ provided that our base field $k$ is finite;  for $k$
finite, the cohomology  of finite functors $F \in \F$ is equal to
the stabilized cohomology of general linear groups with coefficients
determined by $F$ \cite{B2},\cite[App]{FFSS}.

On the other hand, many natural coefficients modules for the general
linear group are not given by functors but by bifunctors
(contravariant in one variable, covariant in the other variable).
Initially motivated by the quest to determine the group cohomology
$\coh^*(\GL(n,\Z/p^2),k)$, efforts have been made to compute the
cohomology of $\GL(n,k)$ with coefficients in symmetric and exterior
powers of the adjoint representation $\gl_n$ (cf. \cite{EF}). These
coefficients are not given by functors but by bifunctors.  In this
paper, we provide computational tools and first computations towards
the determination of the stable (with respect to $n$) values of
$\coh^*(\GL(n,k),S^d(\gl))$ and $\coh^*(\GL(n,k),\Lambda^d(\gl))$.

    Our first task is to formulate in terms of $\Ext$ groups in the category
$\B$ of strict polynomial bifunctors the stable version of rational
cohomology of the algebraic group $\GL$ with coefficients determined
by the given bifunctor. In Theorem \ref{rat}, we show that
\[ \Ext_\B^*(\Gamma^d(\gl),T) ~ \cong \coh^*_{rat}(\GL_n,T(k^n,k^n)) \]
where $T$ is a strict polynomial bifunctor of homogeneous bidegree
$(d,d)$ with $n \geq d$. In the special case that $T$ is of the form
$A(\gl)$ (for example, $S^d(\gl))$, we write this as
\[ \coh^*_\Po(\GL,A) ~ \cong \coh^*_{rat}(\GL_n,T(k^n,k^n)). \]
As for rational cohomology, the most relevant coefficients are given
by beginning with a strict polynomial bifunctor $T$ and applying
the Frobenius twist operation (i.e., $I^{(1)} \circ (-)$)
sufficiently often until the  $\Ext$-group of interest stabilizes.
This ``generic" strict polynomial bifunctor cohomology is our main
target of computations.

    In \cite{FS}, the fundamental computation of
$\Ext_\Po^*(I^{(r)},I^{(r)})$ is achieved, modeled on the
computation of $\Ext_\F^*(I,I)$ in \cite{FLS}. For bifunctor
cohomology, the computation of
\begin{equation}
%\label{equiv}
\coh^*_\Po(\GL,\otimes^{n(r)}) \equiv
\Ext_\B(\Gamma^{np^r},\hc(I^{(r)},I^{(r)})^{\otimes n})
\end{equation}
given in Theorem \ref{computeI} plays an analogous role.

    We prove various useful formal results concerning bifunctor cohomology.
For example, in \S 2 we relate the strict polynomial bifunctor
cohomology with coefficients in a functor of separable type (i.e.,
of the from $\hc(A,B)$ where $A,B$ are strict polynomial functors)
to $\Ext$ computations in the category $\Po$.  In \S 3, we establish
a base change result (one of the important advantages of strict
polynomial functors/bifunctors in contrast to ``usual"
functors/bifunctors) and a twist stability theorem; both results
follow from analogous results proved in \cite{FFSS} for
$\Ext$-groups in the category $\Po$.

    In \S 4, we consider bifunctor cohomology such as
$\coh^*_\Po(\GL,S^d(\gl^{(r)}))$ where $d$ is less than $p$.  This
is in principle completely computable thanks to (\ref{equiv}).
However, computations for $p \leq d$ would appear to be much more
difficult (recently however, Antoine Touz\'e \cite{T} computed
$\coh^*_\Po(\GL,S^d(\gl^{(r)}))$ for all degrees $d$). In \S 5, we
work out the case $p = 2 = d$, a computation which is simplified by
the use of Touz\'e's argument. The applicability of the computation
of \S 5 is extended in \S 6.

    In the remaining two sections, we relate our computations of strict
polynomial bifunctor cohomology to the cohomology of the finite
groups $\GL(n,k)$ where $k$ is a finite field of characteristic $p$.
In \S 7, we develop sufficient formalism for bifunctor cohomology to
enable comparison of this bifunctor cohomology with both the
cohomology of strict polynomial bifunctors and with group
cohomology.  Many explicit computations of group cohomology are
presented in \S 8.

We are most grateful to Antoine Touz\'e who pointed out an error
in our computation of degree 2 bifunctor cohomology in characteristic
2 and provided a simple method to circumvent our difficulty (cf. \S 5.3).
    The first author gratefully acknowledges the hospitality of Northwestern
University where this project was initiated.  Both authors thank
l'Institut Henri Poincar\'e for providing a working environment
which facilitated the completion of this paper.

%%%%%%%%%%%%%%%%%%%%%%%%%%%%%%%
    %%Section1%%%
%%%%%%%%%%%%%%%%%%%%%%%%%%%%%%%%

\section{Rational cohomology of GL as Ext of strict polynomial bifunctors}
\bigskip

We fix a prime $p$ and a field $k$ of characteristic $p$.  Let $\V$
denote the category of finite dimensional $k$-vector spaces and
$k$-linear maps.   In this first section, we introduce the category
of strict polynomial bifunctors and relate $\Ext $-groups in this
category with the rational cohomology of the (infinite) general
linear group.   The reader should be aware that our bifunctors are
contravariant in the first variable, covariant in the second;
bifunctors covariant in each variable were considered for example in
\cite{SFB1}.

In \cite{FS}, the concept of a strict polynomial functor $T: \V \to
\V$ was introduced, a modification of the usual notion of a functor
from $T: \V \to \V$, formulated so that the action of $\GL(W)$ on
$T(W)$ is a rational representation for each $W $ in $ \V$.  Thus,
$T$ consists of the data of an association  of $T(W) \in \V$ for
each $W $ in $ \V$, and a polynomial mapping $\Hom_k(V,W) \to
\Hom_k(T(V),T(W))$ for each pair $V,W $ in $ \V$. [By definition, a
polynomial mapping $V \to W$ between finite dimensional $k$-vector
spaces is an element in $S^*(V^\#) \otimes W$, an element in the
tensor product of the symmetric algebra on the $k$-linear dual of
$V$ with  $W$.] A strict polynomial functor is said to be
homogeneous of degree $d$ if for all pairs $V,W $ in $ \V$ the
polynomial mapping   $\Hom_k(V,W) \to \Hom_k(T(V),T(W)$ is of degree
$d$ (i.e., is an element of $S^d(\Hom_k(V,W)^\#) \otimes
\Hom_k(T(V),T(W)$).

 The abelian category of strict polynomial functors of bounded
degree, denoted $\Po$, is a direct sum of subcategories $\Po_d$ of strict
polynomial functors homogeneous of degree $d$, $\Po = \oplus_d \Po_d$.

We recall a few examples of strict polynomial functors:  the functors
\[
\otimes^d: \V \to V, \quad S^d: \V \to \V, \quad \Gamma^d: \V \to \V, \quad
\Lambda^d: \V \to \V
\]
are each strict polynomial functors homogeneous of degree $d$. Here,
the $d$-th symmetric power $S^d(V)$ of $V\in \V$ is the vector space
of coinvariants of the symmetric group $\Sym_d$ acting on
$\otimes^d(V) = V^{\otimes d}$ by permuting the factors, whereas the
$d$-th divided power $\Gamma^d(V)$ is the vector space of invariants
of $\Sym_d$ acting on $\otimes^d(V)$. For any $W $ in $ \V$ and $T $
in $ \Po_d$, as shown in \cite[2.10]{FS}, the Yoneda lemma yields
natural identifications
\begin{equation}
\label{ident-0}
 \Hom_{\Po_d}(\Gamma^d(\Hom_{k}(W,-),T)~ \cong ~ T(W) ~ \cong ~
\Hom_{\Po_d}(T,S^d(\Hom_k(W,-))^\#.
\end{equation}
Thus, $P_W, I_W$, defined by
\[  P_W ~ =: ~\Gamma^d(\Hom_k(W,-)), \quad I_W =: ~ S^d(\Hom_k(W,-)), \]
are respectively projective and injective objects of $\Po_d$.
Moreover, $P_W$ is a projective generator of $\Po_d$ provided that
$\dim_kW \geq d$, since the natural map
\[ T(W) \otimes \Gamma^d(\Hom_k(W,-),-) \to T \]
is surjective if $\dim_kW \geq d$. Another important strict
polynomial functor, homogeneous of degree $p^r$, is the functor
\[ I^{(r)}: \V \to \V \]
with the property that the structure polynomial maps
are identified with the $p^r$-th power polynomial map.

\begin{defn}
We consider the category $\Po^{op} \times \Po$ of {\it strict polynomial bifunctors}
of bounded degree (contravariant in the first variable, covariant in the second).

Thus, a strict polynomial bifunctor $T$ is a pair of functions, the first of which
assigns to each pair $V,W$ in $\V$ some $T(V,W) \in \V$ and the second of
which assigns to each $V,W,V',W'$ a polynomial map
\begin{equation}
\label{map}
\Hom_k(V,V^\prime) \times \Hom_k(W,W^\prime) \to
\Hom_k(T(V^\prime,W),T(V,W^\prime))
\end{equation}
satisfying the condition that $T(V,-), T(-,W^\#)^\#: \V \to \V$ are
strict polynomial functors for each $V,W $ in $ \V$ of uniformly
bounded degree. If $T$ is in $\Po_d^{op} \times \Po_e$, then we say
that $T$ is homogeneous of bidegree $(d,e)$.
\end{defn}

To emphasize the functorial nature of $T \in \B$, we shall use at times
alternate notations
\[ T~ \equiv ~ T(-,-) ~ \equiv T(-_1,-_2). \]
Observe that if $T$ is a strict polynomial bifunctor, then for any
$W $ in $ \V$ the action
\begin{equation}
\label{action} \GL(W) \times T(W,W) ~ \to ~ T(W,W), \quad (g,x)
\mapsto T(g^{-1},g)(x)
\end{equation}
is {\it rational}.  Namely, for any commutative $k$-algebra $A$,
this action extends to an action $\GL(A\otimes_kW) \otimes
A\otimes_k T(W,W) \to A\otimes_k T(W,W).$

If $T$ is in $\B$ and if $P$ is in $\Po$, then the composite $P
\circ T$ is once again a strict polynomial bifunctor.  If $T_1, T_2$
are strict polynomial functors and $F$ a strict polynomial
bifunctor, then the composite $F(T_1(-),T_2(-))$ is once again a
strict polynomial bifunctor.

For two strict polynomial functors $A_1$ and $A_2$, let $\hc
(A_1,A_2)$ denote the strict polynomial bifunctor defined by
\[ \hc (A_1,A_2)(V,W) = \Hom_k(A_1(V),A_2(W)), \quad V,W \in \V. \]
We refer to such functors as functors of
\emph{separable type}.  One can readily verify the natural identification
\begin{equation}
\label{hom}
\Hom_\B(\hc(A_1,A_2),\hc(B_1,B_2)) ~\cong ~ \Hom_\Po(B_1,A_1)
\otimes \Hom_\Po(A_2,B_2)
\end{equation}

\begin{prop}
\label{enuf}
The category of strict polynomial bifunctors of bounded degree, $\B$ admits
a decomposition $\B \simeq \bigoplus_{d,e} (\Po_d)^{op} \times \Po_e$;
in other words, every
strict polynomial bifunctor of bounded degree can be written naturally as a
direct sum of strict polynomial bifunctors of homogeneous bidegree.

Any $T\in \B$ homogeneous of bidegree $(d,e)$
admits a projective resolution by
(projective strict polynomial) bifunctors of the form
\newline
\[ P_{V,W}^{d,e}(-_1,-_2) = \hc(I_V(-_1),P_W(-_2) ) =
\Gamma^d\Hom_k(-_1,V) \boxtimes \Gamma^e\Hom_k(W,-_2) \] and an
injective resolution by (injective strict polynomial) bifunctors of
the form
\[ I_{V,W}^{d,e}(-_1,-_2) =  \hc(P_V(-_1),I_W(-_2))
 = S^d\Hom_k(-_1,V) \boxtimes S^e\Hom_k(W,-_2). \]
 Here $\boxtimes$ is the external tensor product producing a bifunctor
 from a pair of functors.
\end{prop}

\begin{proof}
Direct sum decomposition of $\B$ is proved exactly as decomposition of
$\Po$ is proved in \cite{FS}.  Using \ref{ident-0}, we easily obtain
the isomorphism
\begin{equation}
\label{respect} \Hom_\B(P_{V,W}^{d,e},T)  \simeq T(V,W) \simeq
\Hom_\B(T,I_{V,W}^{d,e})^\#
\end{equation}
natural with respect to $T \in (\Po_d)^{op}\times \Po_e$. In
particular, each $P_{V,W}^{d,e}$ is projective and each
$I_{V,W}^{d,e}$ is injective in $(\Po_d)^{op}\times \Po_e$. As
argued in \cite[2.10]{FS}, for $T $ in $ (\Po_d)^{op}\times \Po_e$
there is a natural surjective map
\[ P_{V,W}^{d,e} \otimes T(V,W) \to T \]
and a natural injective map
\[ T \to I_{V,W}^{d,e} \otimes T(V,W)^\# \]
whenever $\dim_kV \geq d, \dim_kW \geq e$.
\end{proof}

The strict polynomial bifunctor (of bidegree $(1,1)$)
\[ \gl ~ \equiv ~ \hc(I,I) ~ \equiv ~ \hc(-,-) ~ \in \B \]
plays a special role as we first see in the following proposition.

\begin{prop}\label{gl}
Let $T$ be a strict polynomial bifunctor homogeneous of
bidegree $(d,d)$ and $W \in \V$ satisfy $\dim_k(W) \geq d$. Then
there is a natural identification
\begin{equation}
\label{gl1} \Hom_\B(\Gamma^d\gl,T) ~ \cong ~ \coh^0(\GL(W),T(W,W)),
\end{equation}
where $T(W,W)$ is given the rational $\GL(W)$-module structure of
(\ref{action}).

Moreover, for any $A_1,A_2 $ in $ \Po_d$, there is a natural
identification
\begin{equation}
\label{gl2} \Hom_\B(\Gamma^d\gl,\hc (A_1,A_2)) ~ \cong ~
\Hom_\Po(A_1,A_2).
\end{equation}
\end{prop}

\begin{proof}
For each $W$ in $\V$, one defines a natural transformation:
\begin{equation}\label{natural}
\Hom_\B(\Gamma^d\gl,T) \to \coh^0(\GL(W),T(W,W))
 \end{equation}
which sends a natural transformation
\[
\Phi : \Gamma^d\gl\to T
 \]
to $\Phi(\text{id}_W^{\otimes d})$ in $T(W,W)$.

 The identifications
\[ \begin{aligned}
\Hom_\B(\Gamma^d\gl,I_{U,V}^{d,d}) \cong ~ & \Gamma^d\Hom_k(U,V)^\#
\cong \\
\cong ~ & S^d\Hom_k(V,U)  \cong
\Hom_\Po(\Gamma^d\Hom_k(U,-),S^d\Hom_k(V,-))
\end{aligned} \]
follow immediately from (\ref{respect}) and (\ref{ident-0}); in
particular, this establishes the validity of (\ref{gl2}) for
bifunctors $T$ of the form $I_{U,V}^{d,d}$. For any $W $ in $ \V$,
composing with the evaluation map to
\[
\Hom_{\GL(W)}(\Gamma^d\Hom_k(U,W),S^d\Hom_k(V,W))\cong
\coh^0(\GL(W),I_{U,V}^{d,d}(W,W)).
\]
recovers the above (\ref{natural}) for $T=I_{U,V}^{d,d}$.
 By \cite[3.13]{FS},
this is an isomorphism when $\dim_kW \geq d$. This establishes
(\ref{gl1}) for bifunctors $T$ of the form $I_{U,V}^{d,d}$.

For a general strict polynomial bifunctor $T$ homogeneous of
bidegree $(d,d)$, we apply Proposition \ref{enuf} to obtain an
injective resolution of $T$ by injectives of the form
$I_{U,V}^{d,d}$.   Thus, (\ref{gl1}) follows by the left exactness
of $\Hom_{\B}(\Gamma^d\gl, -)$ and $\coh^0(\GL(W),-)$ together with
the above verification for functors of the form $I_{U,V}^{d,d}$.
Applying (\ref{gl1}) and \cite[3.13]{FS} once again, we conclude
(\ref{gl2}) for $T$ in full generality.
\end{proof}

The existence of enough injectives (and/or projectives) in $\B$
enables us to define $\Ext $-groups in the evident manner.

\begin{defn}
%\label{ext}
For strict polynomial bifunctors $T_1, ~ T_2$ of bounded degree, we
define the $\Ext $-groups $\Ext_\B^*(T_1,T_2 )$ as the derived
functors of $\Hom_\B(T_1,-)$ applied to $T_2$ (or, the derived
functors of $\Hom_\B(-,T_2)$ applied to $T_1$),
\[ \Ext_\B^i(T_1,T_2 )~ = ~ \bR^i\Hom_\B(T_1,T_2). \]

If $A$ is a strict polynomial functor of degree $d$, then we employ
for convenience the following (somewhat misleading) notation:
\begin{equation}
\label{note} \coh_\Po^*(\GL,A) ~ := \Ext_{\B}^*(\Gamma^d\gl,A \circ
\gl),
\end{equation}
\end{defn}

The following theorem relates $\Ext $-groups for bifunctors to the
rational cohomology of the general linear group, thereby extending
Proposition \ref{gl} to positive cohomological degrees.

\begin{thm}
\label{rat} Let $T$ be a strict polynomial bifunctor homogenous of
bidegree $(d,d)$.  If $\dim(W) \geq d$, then there is a natural
isomorphism
\begin{equation}
\label{GL} \Ext_\B^*(\Gamma^d\gl,T) \stackrel {\sim}{\to}
\coh^*_{rat}(\GL(W),T(W,W)),
\end{equation}
where $\coh^i(\GL(W),T(W,W))$ denotes the $i$-th rational cohomology
group of the algebraic group $\GL(W)$ with coefficients in the
rational $\GL(W)$-module $T(W,W)$.

Furthermore, if $F = \hc (A_1,A_2)$ is of separable type
(with $A_1,A_2$ strict
polynomial functors of degree $d$), then there is a natural isomorphism
\begin{equation}
\label{sep} \Ext_\B^*(\Gamma^d\gl,\hc (A_1,A_2)) \simeq
\Ext_\Po^*(A_1,A_2).
\end{equation}
\end{thm}

\begin{proof}
Since both $\coh^*(\GL(W),-(W,W))$ and $\Ext_\B^*((\Gamma^d\gl,-)$
are cohomological $\delta$-functors on $\B$ and since they agree in
degree 0 by Proposition \ref{gl}, to prove (\ref{GL}) it suffices to
verify that
\begin{equation}
\label{van} \coh^i(\GL(W),T(W,W)) = 0, \quad i > 0
\end{equation}
for $T$ of the form $I_{U,V}^{d,d}$.

By \cite[3.13]{FS},
\[ \coh^*_{rat}(\GL(W),I_{U,V}^{d,d}(W,W)) \simeq
\Ext_\Po^i(\Gamma^d\Hom_k(U,-),S^d\Hom_k(V,-)) \] whenever
$\dim_k(W) \geq d$. Since $\Gamma^d(\Hom_k(U,-)$ is projective in
$\Po_d$ (and $S^d\Hom_k(V,-)$ is injective), we conclude these
groups vanish in positive degrees.  This establishes (\ref{GL}).

Now, (\ref{sep}) follows immediately from (\ref{GL}) and \cite[3.13]{FS}
applied to $T = \hc (A_1,A_2)$ .
\end{proof}

    In Proposition \ref{spec}, we give another construction of the
natural isomorphism (\ref{sep}).\par
\vskip .1in

    We shall frequently use the following computation, a fundamental result
of \cite{FS}.

\begin{thm}\cite{FS}
\label{fs}
For any $r \geq 0$, the graded algebra with unit
\[ E_r ~ := ~ \Ext _\Po^*(I^{(r)},I^{(r)}) \]
is the commutative $k$-algebra generated by classes
\[ e_i \in \Ext^{2p^{i-1}}_\Po(I^{(r)},I^{(r)}), \quad 1\leq i \leq r \]
subject to the relations $e_i^p = 0$.
\end{thm}

We have the following adjunction isomorphism,
proved exactly as in  \cite[1.7.1]{FFSS}.

\begin{prop}
\label{adjunct}
Consider the ``diagonal" functor
\[ - \circ D: (\Po^{op})^{ \times m} \times (\Po)^{\times n} ~ \to ~
(\Po^{op})^{ \times m} \times \Po \] sending a strict polynomial
multi-functor $T$ to the functor $T\circ D$ whose value on
$(W_1,\ldots,W_m,V)$ equals $T(W_1,\ldots,W_m,V,\ldots,V)$; also
consider the ``sum" functor
\[ - \circ \bigoplus: (\Po^{op})^{ \times m} \times \Po ~ \to ~
(\Po^{op})^{ \times m} \times (\Po)^{\times n} \] sending a strict
polynomial multi-functor $S$ to the functor $S \circ \bigoplus$
whose value on $(W_1,\ldots,W_m,V_1,\ldots,V_n)$ equals
$T(W_1,\ldots,W_m,\oplus_iV_i)$.  Then $-\circ D,~ - \circ
\bigoplus$ are both exact, and $-\circ D$ is both left and right
adjoint to $-\circ \bigoplus$. Consequently, there are natural
identifications
\[ \Ext_{(\Po^{op})^{ \times m} \times \Po}^*(S,T \circ D) ~\cong ~
\Ext_{(\Po^{op})^{ \times m} \times (\Po)^{\times n}}^*(S \circ
\bigoplus,T) \]
\[ \Ext_{(\Po^{op})^{ \times m} \times \Po}^*(T\circ D,S) ~\cong ~
\Ext_{(\Po^{op})^{ \times m} \times (\Po)^{\times n}}^*( T,S\circ
\bigoplus). \]

Moreover, similar statements apply with $-\circ D$ and $- \circ
\bigoplus$ applied to the contravariant variables of strict
polynomial multi-functors.
\end{prop}

    Theorem \ref{fs} provides the basic ingredient in the computation
of numerous functor cohomology groups (e.g.\cite{FFSS}).   The
following theorem is a first application of this theorem to
bifunctor cohomology. Although this result has probably been
``known" to experts, we know of no written proof (cf. \cite{B}). The
underlying principle of such a computation is to manipulate the
bifunctors involved so that the $\Ext $-computations reduce to
computations of $\Ext $-groups between external tensor products.

\begin{thm}
\label{computeI}
For any $n \geq 1$, we have a $\Sym_n$-equivariant
isomorphism
\begin{equation}
\label{equiv} \coh_\Po^*(\GL,\otimes^{n(r)}) \simeq E_r^{\otimes n}
\otimes k\Sym_n,
\end{equation}
where the action of $\Sym_n$ on the right hand side is by permutation of
the factors of $E_r^{\otimes n} = (\Ext_\Po^*(I^{(r)},I^{(r)}))^{\otimes n}$
and by conjugation on $k\Sym_n$.

\end{thm}

\begin{proof}
We identify $\gl^{\otimes n}$ with $\hc(\otimes^n,\otimes^n)$.  In
this identification, the left action of a permutation $\sigma$ on
the left-hand side is taken to the composition $\sigma^{-1}\circ (-)
\circ\sigma$.

We apply Proposition \ref{rat}  to conclude the isomorphism
\[
\coh_\Po^*(\GL,\otimes^{n(r)}) \simeq
\Ext_\Po^*(\otimes^{n(r)},\otimes^{n(r)}) ~ = ~
\Ext_\Po^*(\otimes^{n(r)},(I^{(r)} \boxtimes \cdots \boxtimes
I^{(r)}) \circ D). \] By Proposition \ref{adjunct}, the right hand
side is isomorphic to
\begin{equation}
\label{kun}
\Ext_{\Po\times \cdots \times \Po}^*(\otimes^{n(r)} \circ \bigoplus,I^{(r)}
\boxtimes \cdots \boxtimes I^{(r)}).
\end{equation}

Expand $(V_1 \oplus \cdots \oplus V_n)^{\boxtimes n} = (\otimes^n
\circ \bigoplus)(V_1,\ldots,V_n)$ as a direct sum of tensor products
and observe that the only summands which are of degree 1 in each
position are the $n!$ terms of the form $V_{\sigma(1)} \boxtimes
\cdots \boxtimes V_{\sigma(n)}$ for some $\sigma $ in $ \Sym_n$.
Thus, the K\"unneth theorem
\[
\Ext_{\Po \times \Po}^*(P_1 \boxtimes P_2,Q_1 \boxtimes Q_2) ~ \simeq ~
\Ext_\Po^*(P_1,Q_1) \otimes \Ext _\Po^*(P_2,Q_2).
 \]
implies (\ref{equiv}) as an additive isomorphism.

The left permutation action on
$\otimes^{n(r)} \circ \bigoplus$ of (\ref{kun})
yields a right action which
permutes the summands of this additive decomposition. We
identify (\ref{equiv}) as a right $\Sym_n$-module as follows.
For each
permutation $\sigma$, we translate the summand
\[
\Ext_{\Po\times \cdots \times \Po}^*(I^{(r)}_{\sigma(1)}\boxtimes
\cdots\boxtimes
I^{(r)}_{\sigma(n)},I^{(r)}
\boxtimes \cdots \boxtimes I^{(r)}).
 \]
indexed by $\sigma$
back to the summand indexed by the identity using the isomorphism
$\sigma^*$; in other words, we identify $x\otimes\sigma$ with
$\sigma^*(x)\otimes 1$
 In this way, the right action is
given by:
\[ \tau^*(x\otimes\sigma)=\tau^*\circ\sigma^*(x)=x\otimes\sigma\tau \]
which is simply the right action of $\Sym_n$ on $k\Sym_n$ in (\ref{equiv}).

We now determine the left action with respect to this
identification. It is enough to consider an elementary tensor
$x=x_1\otimes\dots\otimes x_n $ in $ E_r^{\otimes n} \otimes 1$
representing a class in cohomological degree $s$. Write each $x_i$
as a Yoneda extension in $\Po$
\[
I^{(r)}=Q^0_i\to Q^1_i\to\cdots\to Q^{s_i}_i\to I^{(r)},
 \]
and choose the tensor product of these to represent $x$ as a Yoneda
extension:
\[
I^{(r)} \boxtimes \cdots \boxtimes I^{(r)}=Q^0\to Q^1\to\cdots\to
Q^s\to I^{(r)} \boxtimes \cdots \boxtimes I^{(r)},
 \]
so that $Q^j$ is the sum
\[ \bigoplus_{j_1+\dots
+j_n=j}Q^{j_1}_1\boxtimes\cdots\boxtimes Q^{j_n}_n. \]Let $\sigma$
be a permutation in $\Sym_n$ and let us represent
$x_{\sigma^{-1}(1)}\otimes\dots\otimes x_{\sigma^{-1}(n)}$ by the
tensor of the corresponding Yoneda extensions:
\[
I^{(r)} \boxtimes \cdots \boxtimes I^{(r)}\to Q^1_\sigma\to\cdots\to
Q^s_\sigma  \to I^{(r)} \boxtimes \cdots \boxtimes I^{(r)}.
 \]
By Proposition \ref{adjunct}, the corresponding class in
$\Ext_\Po^*(\otimes^{n(r)},\otimes^{n(r)})$ is obtained by
precomposing with the diagonal $D$. The commutative diagram:
\[ \xymatrix{
    I^{(r)} \otimes \cdots \otimes I^{(r)}\ar[r]\ar[d]^\sigma&
Q^1\circ D\ar[r]\ar[d]^\sigma&\cdots\ar[r]&
Q^i\circ D \ar[r]\ar[d]^\sigma& I^{(r)} \otimes \cdots\otimes
I^{(r)}\ar[d]^\sigma \\
    I^{(r)} \otimes \cdots \otimes I^{(r)}\ar[r]&
Q_\sigma^1\circ D\ar[r]&\cdots\ar[r]&
    Q^i_\sigma\circ D\ar[r]& I^{(r)} \otimes \cdots \otimes I^{(r)} .
                } \]
implies the relation:
\[ \sigma_*(x_1\otimes\dots\otimes
x_n)= \sigma^*(x_{\sigma^{-1}(1)}\otimes\dots\otimes
x_{\sigma^{-1}(n)})
    =x_{\sigma^{-1}(1)}\otimes\dots\otimes
x_{\sigma^{-1}(n)}\otimes\sigma.
 \]

Thus, combining the above left action on the distinguished summand
with the right $\Sym_n$-action permuting the summands
gives the asserted action.
\end{proof}

%%%%%%%%%%%%%%%%%%%%%%%%%%%%%
% Please feel free to modify!

As we see in the following proposition, Theorem 1.8 together with an
analysis of the $\Sym_3$-module
\[
\coh_\Po^*(\GL,\otimes^{3(r)}) \simeq E_r^{\otimes 3} \otimes
k\Sym_3.
 \]
gives us computations of certain bifunctor cohomology associated to
$\otimes^{3(r)}$. In section \ref{schur}, we apply Theorem
\ref{computeI} and compute the bifunctor cohomology of summands of
$\otimes^{n(r)}$ for any $n \geq 2$.

%
%%%%%%%%%%%%

\begin{prop}\label{ps}
We consider a field $k$ of characteristic $p > 0$.
\begin{itemize}
\item  The Poincar\'e series for $\coh_\Po^*(\GL,\otimes^{3(r)})$ equals
\[ 6\left(  \frac{1-t^{2p^r}}{1-t^2}\right)^3. \]
\item  If $p > 3$, the Poincar\'e series for
$\coh_\Po^*(\GL,\Lambda^{(3(r)})$ equals
\[
\left(  \frac{1-t^{2p^r}}{1-t^2}\right)^3 - 2\left(
\frac{1-t^{4p^r}}{1-t^4}\right) \left(
\frac{1-t^{2p^r}}{1-t^2}\right) + 2\left(
\frac{1-t^{6p^r}}{1-t^6}\right).
 \]
\item  If $p > 3$, the Poincar\'e series for $\coh_\Po^*(\GL,S^{(3(r)})$
equals
\[
2\left(  \frac{1-t^{4p^r}}{1-t^4}\right) \left(
\frac{1-t^{2p^r}}{1-t^2}\right) + \left(
\frac{1-t^{6p^r}}{1-t^6}\right).
 \]
\end{itemize}
\end{prop}

\begin{proof}
We immediately verify that the Poincar\'e series for $E_r$ equals
$P_r(t) = \frac{1-t^{2p^r}}{1-t^2}$.

Index the isomorphism classes of transitive $\Sym_3$-sets by their
cardinality: $S_1$, the trivial $\Sym_3$-set, $S_2$, $S_3$,
$S_6=\Sym_3$. Every $\Sym_3$-set is uniquely isomorphic to a
disjoint union of these. Finite disjoint unions of copies of these
four sets form a commutative semi-ring under disjoint union and
cartesian product, with $S_1$ as unit. The product is given
explicitly by:
\[
S_1\times S_i=S_i ,\quad S_2\times S_3=S_6,\quad S_i\times S_i = i
S_i, \quad S_i\times S_6=iS_6.
 \]
For example, with respect to the conjugation action, $\Sym_3\simeq
S_1\coprod S_2\coprod S_3$ corresponding to the partition of
$\Sym_3$ into conjugacy classes. Similarly, every finite-dimensional
permutation module is isomorphic to a finite direct sum of $k[S_i]$,
and every graded permutation module which is finite dimensional in
each degree has a Poincar\'e series with coefficients in the above
semi-ring.

For our purpose, the relevant example is the third tensor power,
$E_r^{\otimes 3}$, of the graded vector space $E_r$ with Poincar\'e
series $P_r(t)=\sum_i\dim E_r^i t^i$.  A direct computation verifies
that the Poincar\'e series of $E_r^{\otimes 3}$ as a $\Sym_3$-module
equals
\[
S_1P_r(t^3)+S_3(P_r(t^2)P_r(t)-P_r(t^3))+
\frac{S_6}{6}(P_r(t)^3-3P_r(t^2)P_r(t)+2P_r(t^3)).
 \]
(The subspace fixed by $\Sym_3$ is spanned by elements of the form
$x\otimes x \otimes x $ in $ E_r^{\otimes 3}$, leading to the first
summand; the cycles of length 3 are represented by elements of the
form $x \otimes x \otimes y, ~ x \not= y$, which leads to the second
summand; the cycles of length 6 are represented by elements of the
form $x\otimes y \otimes z$ with $x,y,z$ distinct.) After
multiplication using the product given explicitly above, one gets
that the Poincar\'e series for $E_r^{\otimes 3} \otimes k\Sym_3$ is
equal to
\[
S_1P_r(t^3)+S_2P_r(t^3)+S_3(4P_r(t^2)P_r(t)-3P_r(t^3))+
S_6(P_r(t)^3-2P_r(t^2)P_r(t)+P_r(t^3)).
 \]

The Poincar\'e series of the graded vector spaces obtained by
applying symmetrization and antisymmetrization functors (the
functors $s^3(-)$) and $\lambda^3(-)$ of \S \ref{schur} below) to
$E_r^{\otimes 3} \otimes k\Sym_3$ is then obtained term by term from
the corresponding result for the permutation modules $k[S_i]$. If
$p$ is not equal to $3$, then
$\lambda^3(k[S_1])=\lambda^3(k[S_3])=0$, $\lambda^3(k[S_2]) =
\lambda^3(k[S_6]) = k$, and $s^3(k[S_i]) = k, ~ 1 \leq i \leq 6$.

As a result, the Poincar\'e series of the graded vector space
$\lambda^3(E_r^{\otimes 3} \otimes k\Sym_3)$ is equal to
\[
P_r(t)^3-2P_r(t^2)P_r(t)+2P_r(t^3) \] whereas the Poincar\'e series
for $s^3(E_r^{\otimes 3} \otimes k\Sym_3)$ is equal to
\[  2P_r(t^2)P_r(t)+ P_r(t^3). \]
By Section \ref{schur} below, this yields the desired result.
\end{proof}

%%%%%%%%%%%%%%%%%%%%%%%%%%%%%%%%%%%%%
%Section 2
%%%%%%%%%%%%%%%%%%%%%%%%%%%%%%%%%%%%%%

\section{Cohomology of bifunctors of separable type}

    In Theorem \ref{rat}, we established the isomorphism (\ref{sep})
between $\Ext $-groups in the category $\B$ for bifunctors of
separable type and $\Ext $-groups in the category $\Po$ of strict
polynomial functors.  The purpose of this section is to further
exploit the special properties of such bifunctors of separable type.
% We are especially interested in Yoneda products. In particular, the
% classical Lemma \ref{Hopf} states that a certain precomposition with
% $gl$ is the coproduct of a graded Hopf algebra structure, a
% verification which will play a central role in the calculation of \S
% \ref{degree2}.

We begin with the following elementary lemma.

%%%%%%%%%%%%
%% I straighten out the lemma to match the spectral sequence below
%% Anyway, it was not correct
%% the variance in T doesn't match
%% didn't VdK tell us?
%%  see (1.2.1)=(\ref{ident-0}) as well
%%%%%%%%%%%%%%%%%%%%%%%%

\begin{lemma}
\label{switch} For $V,W $ in $ \V$ and $T $ in $ \B$, there is a
natural identification
%% \[ \Hom_\B(T,\hc(P_V,I_W)) ~ \cong ~
%% \Hom_\Po(\Hom_\Po(T(-_1,-_2),P_V(-_1)),I_W(-_2)). \]
\[ \Hom_\B(T,\hc(P_V,I_W)) ~ \cong ~
\Hom_\Po(P_V(-_1),\Hom_\Po(T(-_1,-_2),I_W(-_2))). \]
\end{lemma}

\begin{proof}
By (\ref{respect}),
%% \[ \Hom_\B(T,\Hom_\Po(P_V,I_W)) ~ \cong ~ T(V,W). \]
\[ \Hom_\B(T,\Hom_\Po(P_V,I_W)) ~ \cong ~ T(V,W)^\#. \]
On the other hand,
%\[ \Hom_\Po(T(-_1,-_2),P_V(-_1)) ~ \cong ~T(V,-_2), \]
\[ \Hom_\Po(T(-_1,-_2),I_W(-_2)) ~ \cong ~T(-_1,W)^\#, \]
by (\ref{ident-0}); thus, another application of (\ref{ident-0})
implies
%% \[ \Hom_\Po(\Hom_\Po(T(-_1,-_2),P_V(-_1)),I_W(-_2)) ~ \cong ~ T(V,W). \]
\[ \Hom_\Po(P_V(-_1),\Hom_\Po(T(-_1,-_2),I_W(-_2)))~ \cong ~ T(V,W)^\#. \]

The composite of these natural identifications is readily seen to be
the identification asserted in the statement of the lemma.
\end{proof}

The preceding lemma leads to the following spectral sequence, which
provides an explanation and an extension of (\ref{sep}).

\begin{prop}
\label{spec}
For two strict polynomial functors $A_1$ and $A_2$ homogeneous of degree $d$,
we consider the bifunctor
$\hc (A_1,A_2)$ of separable type.  For any strict polynomial bifunctor $T$, there is a
convergent spectral sequence of the form
\begin{equation}
\label{spec1}
E_2^{s,t} ~ = ~ \Ext_{\Po}^s(A_1(-_1),\Ext_{\Po}^t(T(-_1,-_2),A_2(-_2)))
~ \Rightarrow ~ \Ext_{\B}^{s+t}(T,\hc (A_1,A_2)),
\end{equation}
natural in $A_1$, $A_2$ and $T$.

If $T = \Gamma^d \gl $, then this spectral sequence collapses to
give the natural isomorphism
\[ \Ext_{\B}^*(\Gamma^d \gl,\hc (A_1,A_2)) \simeq  \Ext_{\Po}^*(A_1,A_2). \]
\end{prop}

\begin{proof}
If $P_\bu \to A_1 ~ \in \Po$ is a projective resolution of $A_1$ and if
$A_2 \to I^\bu ~ \in \Po$ is an injective resolution of $A_2$, then the
double complex $\hc(P_\bu,I^\bu)$ of injective bifunctors in $\B$ has
total complex which is an injective resolution of
$\hc (A_1,A_2)$.  Thus, the cohomology of the total complex of
the bicomplex $\Hom_\B(T,\hc(P_\bu,I^\bu))$ equals
$\Ext_{\B}^*(T,\hc (A_1,A_2))$.  On the other hand, Lemma \ref{switch}
identifies this bicomplex with
\[ \Hom_\Po(P_\bu(-_1),\Hom_\Po(T(-_1,-_2),I^\bu(-_2)). \]
If we take the iterated cohomology of this bicomplex first with
respect to the variable index of $I^\bu$ we obtain
\[ \Hom_\Po(P_\bu(-_1),\Ext^t_\Po(T(-_1,-_2),A_2)); \]
then taking cohomology with respect to the variable index of $P_\bu$,
we obtain the asserted $E_2$-term
\[ \Ext_{\Po}^s(A_1(-_1),\Ext_{\Po}^t(T(-_1,-_2),A_2(-_2))). \]
Thus, (\ref{spec1}) is one of the usual spectral sequences
associated to the bicomplex $\Hom_\B(T,\hc(P_\bu,I^\bu))$.

For $T = \Gamma^d \gl $, (\ref{gl2})  provides the identification of
bicomplexes 
\begin{equation}\label{identification:eq}
 \Hom_\B(
\Gamma^d \gl,\hc(P_\bu,I^\bu)) ~ \cong ~\Hom_\Po(P_\bu,I^\bu);
\end{equation}
the cohomology of its total complex is thus equal to
$\Ext_{\Po}^*(A_1,A_2)$; the constructed spectral sequence clearly
collapses (since $\Ext^i(P_s,I^\bu) = 0, ~ i > 0$).
\end{proof}

    The following proposition, extending (\ref{hom}) to
all cohomology degrees, can be viewed as a K\"unneth Theorem for
bifunctor cohomology.

\begin{prop}
\label{Ku} Let $A_1,A_2,B_1,B_2$ be strict polynomial functors of
bounded degree. Then the spectral sequence (\ref{spec1}) with $T =
\hc(B_1,B_2)$ collapses to yield the natural identification
\[
\Ext_{\B}^*(\hc (B_1,B_2),\hc (A_1,A_2))~ \cong ~
\Ext^*_{\Po}(A_1,B_1)\otimes\Ext_{\Po}^*(B_2,A_2).
 \]
\end{prop}

\begin{proof}
Consider the bicomplex $\Hom_\B(\hc(A_1,A_2),\hc(P_\bu,I^\bu))$;
its total complex computes
$\Ext_{\B}^*(\hc (A_1,A_2),\hc (B_1,B_2)).$
By (\ref{hom}), this bicomplex can be identified
with $\Hom_\Po(P_\bu,B_1)\otimes \Hom_\Po(B_2,I^\bu)$.
The iterated cohomology of this bicomplex equals its total cohomology
as well as equals
$\Ext^*_{\Po}(A_1,B_1)\otimes\Ext_{\Po}^*(B_2,A_2)$
\end{proof}

\begin{rem}\label{products}
The isomorphism of Proposition \ref{Ku} is compatible with Yoneda
products.  In particular, there is a ring isomorphism
\[
\Ext_{\Po}^*(A_1,A_1)^{op}\otimes\Ext_{\Po}^*(A_2,A_2)
\stackrel{\cong}{\longrightarrow} \Ext_{\B}^*(\hc (A_1,A_2),\hc
(A_1,A_2)). \]

The Yoneda product defines a right action of $\Ext_{\Po}^*(A_1,A_1)$
and a left action of $\Ext_{\Po}^*(A_2,A_2)$ on  the $E_2$-term of the
spectral sequence in Proposition \ref{spec}, while
the abutment of the spectral sequence is a module over the ring
\[ \Ext_{\B}^*(\hc (A_1,A_2),\hc (A_1,A_2)). \]
Because of the way the
spectral sequence of Proposition \ref{spec} is constructed,
it is compatible with such Yoneda products,
which makes it a $\Ext_{\Po}^*(A_1,A_1)$-$\Ext_{\Po}^*(A_2,A_2)$
bi-module spectral sequence.
\end{rem}

%% %%%%%%%%%%
%% We just suppress whole of 2.5, 2.6 with proof
%%%%%%%%%%%%%%%

%%%%%%%%%%%%%%%%%%%%%%%%%%%%
%Section 3
%%%%%%%%%%%%%%%%%%%%%%%%%%%%%

\section{Base change and twist stability for bifunctor cohomology}
\subsection{Base change}

    If $K/k$ is a field extension of our chosen base field $k$, then we
let $\V_K$ denote the category of finite dimensional $K$-vector spaces
and $\Po_K$ the category of strict polynomial functors on $\V_K$
of bounded degree.  We denote by $\B_K$ the category of strict polynomial
bifunctors on $\V_K^{op} \times V_K$.

    Following the construction in \cite[2.5]{SFB1}, we define the base change
$T_K \in \B_K$ of a strict polynomial bifunctor $T \in \B$ as
follows.  For any $V $ in $ \V$, let $V_K$ denote $V\otimes_k K$,
the base change of $V$ to $K/k$.  Then we define
\[ T_K(V^\prime,W^\prime) ~ := ~
\varinjlim_{\V_K(V^\prime)^\sim \times \V_K(W^\prime)} T(V,W) \]
where $\V_K(V^\prime)^\sim$ is the category whose objects are pairs
$(V,\phi), ~ V \in \V, ~ \phi: V^\prime \to V\otimes_k K $ and whose
maps $(V,\phi) \to (V_1, \phi_1)$ are $K$-linear maps $\theta: V
_1\otimes_k K \to V \otimes_k K$ such that $\phi =  \theta\circ
\phi_1$; and where $\V_K(W^\prime)$ is the category whose objects
are pairs $(W,\psi), ~ W \in \V, ~ \psi: W\otimes_k K \to W^\prime$
and whose maps $(W,\psi) \to (W_1,\psi_1)$ are $K$-linear maps
$\rho: W\otimes_k K \to W_1 \otimes_k K$ such that $\psi = \rho
\circ \psi_1$. As in \cite[2.5]{SFB1}, $T_K$ is a well-defined
strict polynomial bifunctor whose value on a pair of the form
$(V_K,W_K)$ equals $T(V,W)_K$.

    Cohomological base change for strict polynomial bifunctors is formulated
in the following proposition.

\begin{prop}
Let $K/k$ be a field extension and $S,T \in \B$ be strict polynomial
bifunctors of bounded degree.  Then there is a natural isomorphism
of graded $K$-vector spaces
\[ \Ext_\B^*(S,T) \otimes_k K ~ \stackrel{\sim}{\to} \Ext_{\Po_K^{op} \times
 \Po_K}^*(S_K,T_K). \]
\end{prop}

\begin{proof}
As observed in \cite[2.6]{SFB1}, sending $T$ to $T_K$ is exact. The
verification that this base change preserves projective objects of
the form
\[ \Gamma^d\Hom_k(-_1,W) \boxtimes
\Gamma^d(\Hom_k(V,-_2) \] follows from the fact also verified in
\cite[2.6]{SFB1} that the base change of $\Gamma^d\Hom_k(V,-) \in
\Po$ is naturally isomorphic to $\Gamma^d\Hom_K(V_K,-)$.  Then the
natural map
\[ \hc_\B(S,T) \otimes_k K ~ \stackrel{\sim}{\to}
\hc_{\Po_K^{op} \times \Po_K}(S_K,T_K) \] is an isomorphism for $T$
of the form $\Gamma^d\Hom_k(-_1,W) \boxtimes
\Gamma^d(\Hom_k(V,-_2)$. Since any $T$ admits a resolution by such
projective objects, the proposition follows.
\end{proof}

\subsection{Frobenius twist}
For any strict polynomial bifunctors  $T$, the exactness of $I^{(1)}
\circ (-)$ determines the {\it Frobenius twist} map
\begin{equation}
\label{tw1} \Ext _\B^*(S,T) ~ \to ~ \Ext _\B^*(S^{(1)},T^{(1)}).
\end{equation}

 {\it Twist injectivity} for strict polynomial bifunctor cohomology,
given below, is essentially a reformulation of a theorem for
rational cohomology proved by H. Andersen.

\begin{prop}
\label{inj} For any strict polynomial bifunctor $ T \in \B$ of
bidegree $(d,d)$, the composition of (\ref{tw1}) (taking $S =
\Gamma^d\gl$) and the natural map induced by $\Gamma^{pd} \to
\Gamma^{d(1)}$ (dual to the $p$-th power map  $S^{d(1)} \to
S^{pd}$),
\begin{equation}
\label{tw2} \Ext _\B^*(\Gamma^d\gl,T) \to \Ext
_\B^*(\Gamma^{pd}\gl,T^{(1)}).
\end{equation}
is injective.
\end{prop}

\begin{proof}
For any $W $ in $ \V$,  consider the square
\[
\xymatrix{\Hom_\B(\Gamma^d\gl,T) \ar[d]\ar[r]^\cong &
\coh^0(\GL(W),T(W,W)) \ar[d]\\
\Hom_\B(\Gamma^{pd}\gl,T^{(1)}) \ar[r]^\cong &
\coh^0(\GL(W),T(W,W)^{(1)}), }
 \]
where the horizontal isomorphisms are those of (\ref{GL}), the left
vertical map is the map induced by twist and $\Gamma^{pd} \to
\Gamma^{d(1)}$, and the right vertical map is the Frobenius twist
map.  Tracing through the identification given in the proof of
Proposition \ref{gl} for $T = I_{U,V}^{d,d}$, we easily conclude
that this square commutes for such $T$.  For general $T$, we
consider the beginning of a resolution of $T$ by injectives of this
form, $T \to I^0 \to I^1$ and use the left exactness of
$\Hom_B(\Gamma^d\gl,-)$ and $\coh^0(\GL(W),-)$. Since the right
vertical map is evidently injective, so is the left vertical map;
this verifies the asserted injectivity in cohomological degree 0.

For higher cohomological degree, embed $T$ in some injective strict
polynomial bifunctor $J$ and denote the quotient $J/T$ by $\ol J$.
We identify (\ref{tw2}) in cohomological degree 1 with the map on
quotients induced by maps in cohomological degree 0, since
\[
\Ext
_\B^1(\Gamma^d\gl,T) ~ = ~ \Hom_\B(\Gamma^d\gl,\ol J)/
\Hom_\B(\Gamma^d\gl,J) \] and
\[ \coh^1(\GL,T(W,W)) ~ = ~ \coh^0(\GL,\ol J(W,W))/ \coh^0(\GL,J(W,W)).
\]
In cohomological degree $i> 1$, we identify (\ref{tw2}) with the map in
cohomological degree $i-1$ with $T$ replaced by $\ol J$.  Consequently,
we conclude the commutativity of the square
\begin{equation}
\label{ext}
\xymatrix{\Ext^*_\B(\Gamma^d\gl,T) \ar[d]\ar[r]^\cong &
\coh^*(\GL(W),T(W,W)) \ar[d]\\
\Ext^*_\B(\Gamma^{pd}\gl,T^{(1)})) \ar[r]^\cong &
\coh^*(\GL(W),T(W,W)^{(1)}) }
\end{equation}
where the left vertical map is the map of the assertion, the right
vertical map is the twist map in rational cohomology, and where the
horizontal  isomorphisms are those (\ref{GL}).

    The injectivity of the left vertical map thus follows from Andersen's
theorem asserting that the right vertical map is injective (cf.
\cite{A}, \cite[II.10.16]{J}).
\end{proof}

    We now prove {\it twist stability} for strict polynomial bifunctor
cohomology, basically by reducing the assertion to the theorem of
``strong twist stability" given in \cite{FFSS}.

\begin{thm}
\label{twist} For any $S,~T $ in $ \B$ of bidegree $(d,d)$, the
twist map of (\ref{tw1}) in cohomological degree $s$,
\begin{equation}
\label{tw3} \Ext _\B^s(S^{(r)},T^{(r)}) ~ \to ~ \Ext
_\B^s(S^{(r+1)}\gl,T^{(r+1)}),
\end{equation}
is an isomorphism provided that $r \geq log_p(\frac{s+1}{2}).$

Moreover, the twist map of (\ref{tw2}) in cohomological degree $s$,
\begin{equation}
\label{tw4} \Ext _\B^s(\Gamma^{p^rd}\gl,T^{(r)}) \to \Ext
_\B^s(\Gamma^{p^{r+1}d}\gl,T^{(r+1)}),
\end{equation}
is also an isomorphism
provided that $r \geq log_p(\frac{s+1}{2})$.
\end{thm}

\begin{proof}
We first assume that $S = \hc(A_1,A_2),~
T = \hc(B_1,B_2)$ are of separable type and identify (\ref{tw1}) with
\begin{equation}
\label{tw5}
\Ext_\Po^*(B_1,A_1)\otimes \Ext_\Po^*(A_2,B_2) ~ \to ~
\Ext_\Po^*(B_1^{(1)},A_1^{(1)})\otimes \Ext_\Po^*(A_2^{(1)},B_2^{(1)}).
\end{equation}
Provided that $ r \geq log_p(\frac{s+1}{2})$, (\ref{tw5}) is an
isomorphism in cohomological degree $s$ by \cite[4.10]{FFSS}.

More generally, assume that $S = \hc(A_1,A_2)$ is a bifunctor of
separable type but that $T$ is an arbitrary strict polynomial
bifunctor homogeneous of degree $(d,d)$.   Choose a resolution of
$T$, $T\to J^\bu$, by bifunctors of separable type (each of which is
strict polynomial of bidegree $(d,d)$) and compare the map induced
by twist of hyperext spectral sequences obtained by applying
$\Hom_\B(S^{(r)},-)$ and $\Hom_\B(S^{(r+1)},-)$ to double complexes
which are injective resolutions of $J^{\bu (r)}$ and $J^{\bu
(r+1)}$. The preceding special case implies that we have an
isomorphism on $E_1$-terms of these spectral sequences,
\[
E_1^{s,t} ~ \stackrel{\sim}{\to} ~ {}^\prime E_1^{s,t}, \quad r \geq
log_p \frac{t+1}{2}. \] We identify the map on abutments,
\[ \Ext _\B^*(S^{(r)},T^{(r)}) \cong \bE_\B^*(P_\bu^{(r)},J^{\bu (r)}) ~\to
\]
\[
\bE_\B^*(P_\bu^{(r+1)},J^{\bu (r+1)}) \cong \Ext _\B^*(S^{(r+1)},T^{(r+1)}),
\]
in cohomological degree $s$ with (\ref{tw3}), so that an easy spectral
sequence argument completes the proof of the first assertion in this case.

Finally, we allow $S$ to be an arbitrary strict polynomial bifunctor
homogeneous of degree $(d,d)$ and choose a resolution $P_\bu \to S$
by bifunctors of separable type.  Comparing hyperext spectral
sequences as in the previous paragraph enables us to conclude the
validity of the first assertion.

To prove the second assertion, we proceed in a similar fashion.  We
first observe that the commutativity of (\ref{ext}) in the special
case that $T$ is of separable type together with (\ref{sep}) enables
us to apply the ``strong twist stability theorem" of
\cite[4.10]{FFSS} to conclude that (\ref{tw4}) is an isomorphism
whenever $ r \geq log_p \frac{s+1}{2}.$

More generally, we choose a resolution $T \to J^\bu$ of $T$ by
strict polynomial bifunctors of separable type.  We consider the map
of hyperext spectral sequences from
\[ E_1^{s,t} = \Ext_\B^t(\Gamma^{p^rd}\gl,J^{s(r)})
\Rightarrow \Ext_\B^{s+t}(\Gamma^{p^rd}\gl,T^{(r)}) \] to
\[
 {}^\prime E_1^{s,t} = \Ext_\B^t(\Gamma^{p^{r+1}d}\gl,J^{s(r+1)})
 \Rightarrow \Ext_\B^{s+t}(\Gamma^{p^{r+1}d}\gl,T^{(r+1)}).
 \]
Once again,  an easy spectral sequence comparison theorem implies that
\[ E_1^{s,t} ~ \stackrel{\sim}{\to} ~ {}^\prime E_1^{s,t}, \quad r \geq
log_p \frac{t+1}{2}, \] so that the induced map on abutments
$E_\infty^s \to {}^\prime E_\infty^s$ is an isomorphism provided
that $r \geq log_p(\frac{s+1}{2})$.
\end{proof}

We restate Theorem \ref{twist} in the special case in which $T =
A\circ \gl$, using the notation of (\ref{note}).

\begin{cor}
\label{rtwist} For any $A $ in $ \Po$,
\[ \coh^s(\GL,A^{(r)}) ~\to~ \coh^s(\GL,A^{(r+1)}) \]
is an isomorphism provided that $r \geq log_p(\frac{s+1}{2})$.
\end{cor}

%%%%%%%%%%%%%%%%%%%%%%%%%%%%%%
%Section 4
%%%%%%%%%%%%%%%%%%%%%%%%%%%%%%

\section{Cohomology of summands of $gl^{(r)\otimes n}, ~n < p$}
\label{schur}
\bigskip

In this section, we consider summands of the bifunctor
\[ \gl^{\otimes n(r)} ~ := ~ \hc(I,I)^{\otimes n (r)} \]
for $n < p$.  We are particularly interested in the $n$-th symmetric
power $S^n(\gl^{(r)})$ and the $n$-th exterior power
$\Lambda^n(\gl^{(r)})$, but the method of determining the cohomology
of these summands applies more generally to summands given by
idempotents in the group ring $k\Sym_n$.  Our hypothesis $n < p$
implies that this group ring is semi-simple.   We refer the reader
to \cite{FH} for further details.

    Let $\lambda = (\lambda_1,\ldots,\lambda_j)$ be a partition of
$n$ and choose a Young diagram whose shape is given by $\lambda$ (so
that the first row of the Young diagram has $\lambda_1$ ``boxes",
the second has $\lambda_2$ boxes, etc).  Let $t_{\lambda}$ be an
arbitrary choice of a Young tableau (i.e., an assignment of the
numbers $\{ 1,\ldots,n \}$ to these boxes which is strictly
 increasing both to the left and downwards)
with  underlying Young diagram given by $\lambda$. Define $a_\lambda
\in k\Sym_n$ to be the sum (in $k\Sym_n$) of the elements in
$\Sym_n$ which preserve the rows of this tableau, $b_\lambda$ to be
the alternating sum of the elements in $\Sym_n$ which preserve the
columns.  Then some multiple of
\[ c_\lambda := a_\lambda \cdot b_\lambda \in k\Sym_n \]
is idempotent and
\[
\{ V_\lambda := k\Sym_n \cdot c_\lambda; \quad \lambda
~{\text a ~partition~of~} n\} \]
constitutes a complete list of the irreducible representations of $\Sym_n$
(where, as always in this section, we assume $p > n$).

    In particular, if $\lambda = (n)$, then  $V_{(n)} ~ = ~ k$ is the
1-dimensional trivial representation of $\Sym_n$; if $\lambda =
(1,\ldots,1)$, then $V_{(1,\ldots,1)} ~ = ~ k_{sgn}$ is the
1-dimensional sign representation.

For any $k\Sym_n$-module $W$ and any partition $\lambda$ of $n$, we define
\[ s^{\lambda}(W) =: ~ ( k\Sym_n \cdot c_\lambda) \otimes_{k\Sym_n} W. \]
If $W = V^{\otimes n}$ for some finite dimensional vector space $V$ and if
$\Sym_n$ acts on $W$ by permuting the tensor factors, we introduce the
notation
\[ S^\lambda(V) := ~ s^\lambda(V^{\otimes n}); \]
in particular, $S^{(n)}(V) = S^n(V)$ and
$s^{(1,1,\ldots,1)}(V) = \Lambda^n(V)$.

    The following proposition provides an (implicit) description of
cohomology of numerous bifunctors granted the explicit description
of
\[
\coh_\Po^*(\GL,\gl^{\otimes^n(r)}) = E_r^{\otimes n} \otimes k\Sym_n
\]
as a $\Sym_n$-module given in Theorem \ref{computeI}.

\begin{prop}
\label{partition}
Assume $p > n$ and let $\lambda$ be a partition of $n$.  Then
\begin{equation}
\label{sl} \coh_\Po^*(\GL,S^\lambda(\gl^{(r)})) ~ \cong
s^\lambda(\coh_\Po^*(\GL,\gl^{\otimes n(r)})).
\end{equation}
In particular,
\[
\coh_\Po^*(\GL,S^{n(r)}\gl) ~ \cong s^{(n)} (\coh_\Po^*(\GL,\gl^{\otimes n(r)})) ~
\cong ~ k \otimes_{k\Sym_n} \coh_\Po^*(\GL,\gl^{\otimes n(r)}),
\]
\[
\coh_\Po^*(\GL,\Lambda^{n(r)}\gl) ~ \cong
s^{(1,\ldots,1)} (\coh_\Po^*(\GL,\gl^{\otimes n(r)}))
~ \cong ~k_{sgn}\otimes_{k\Sym_n} \coh_\Po^*(\GL,\gl^{\otimes
n(r)}),
\]
\end{prop}

\begin{proof}
The fact that $c_\lambda$ is a quasi-idempotent means that some
non-zero multiple of $c_\lambda$, $e_\lambda = a_\lambda c_\lambda,
$ is idempotent. Thus,
\[
k\Sym_n ~ \cong ~ k\Sym_n\cdot e_\lambda \oplus k\Sym_n \cdot (1-e_\lambda).
\]
By functoriality,
\[ s^\lambda(\coh_\Po^*(\GL,\gl^{\otimes n(r)})) =: ~
\coh_\Po^*(\GL,\gl^{\otimes n(r)}) \cdot c_\lambda ~ = ~
\coh_\Po^*(\GL,\gl^{\otimes n(r)})\cdot e_\lambda \] equals
\[ \coh_\Po^*(\GL,(\gl^{\otimes n(r)})\cdot e_\lambda) ~=~
 \coh_\Po^*(\GL,(\gl^{\otimes n(r)})\cdot c_\lambda) ~=~
 \coh_\Po^*(\GL,S^\lambda(\gl^{(r)})) . \]
\end{proof}

As a simple corollary of Proposition \ref{partition}, we have the
following vanishing of strict polynomial bifunctor cohomology for
$r=0$.

\begin{cor}
Assume $p > n$ and let $\lambda$ be a partition of $n$.  Then
\[ \coh_\Po^s(\GL,S^\lambda(\gl))~ = ~0, \quad s > 0. \]
\end{cor}

\begin{proof}
This follows immediately from Proposition \ref{partition} and the
observation that $E_0$ vanishes in positive degrees.
\end{proof}

We make explicit the case $n=2$ of Proposition \ref{partition}.

\begin{cor}
If $p > 2$, then $c_{(2)} = 1 +\tau , ~ c_{(1,1)} = 1 -\tau  \in k\Sym_2$.    Thus,
\[ \coh_\Po^*(\GL,S^{2(r)}) ~ \cong ~ S^2(E_r) \oplus S^2(E_r),
\quad \coh_\Po^*(\GL,\Lambda^{2(r)}) ~ \cong ~ \Lambda^2(E_r) \oplus
\Lambda^2(E_r). \]

\end{cor}

\begin{proof}
Recall the triviality of the action of $\Sym_2$ on the tensor factor
$k\Sym_2$ of
\[
\coh_\Po^*(\GL,\gl^{\otimes^2(r)}) \cong E_r^{\otimes 2} \otimes k\Sym_2.
\]
Thus, $s_\lambda(\coh_\Po^*(\GL,\otimes^{2(r)}))$ is naturally
isomorphic to two copies of $s_\lambda(E_r^{\otimes 2})$, where the
action of $\Sym_2$ on $E_r^{\otimes 2}$ is permutation of tensor
factors by Theorem \ref{computeI}.
\end{proof}

We derive one more simple corollary which follows immediately from Theorem
\ref{computeI}, Proposition
\ref{partition} and the fact that the graded group $E_r^{\otimes n}$ is
1-dimensional in degree 0.

\begin{cor}
Assume $p > n$ and let $\lambda$ be a partition of $n$.  Then
\[
\coh_\Po^0(\GL,S^\lambda(\gl^{(r)})) ~ \cong s^\lambda(k\Sym_n)
 \]
where $k\Sym_n$ is a $\Sym_n$-module via the conjugation action.
\end{cor}

%%%%%%%%%%%%%%%%%%%%%%%%%%%%%
%Section 5
%%%%%%%%%%%%%%%%%%%%%%%%%%%%

\section{Computations for $S^2(gl^{(r)}), ~ \Lambda^2(gl^{(r)}),
and ~\Gamma^2(gl^{(r)}) ~ \text{for}~p = 2$}
\label{degree2}
\bigskip

In this section, our base field $k$ is an arbitrary field of
characteristic 2  and $r$ a non-negative integer.  This section is
dedicated to verifying the following computation.

\begin{thm}
\label{deg2}
Let $k$ be a field of characteristic 2 and $r \geq 0$ a non-negative integer.
Then
\begin{itemize}
\item The $\Sym_2$-module $\coh_\Po^*(\GL,\otimes^{2{(r)}})$ is isomorphic to
$E_r^{\otimes 2} \otimes k\Sym_2$ with $\Sym_2$-action on
$E_r^{\otimes 2}$ given by permuting the tensor factors and with $\Sym_2$-action
on  $k\Sym_2$ trivial; its Poincare series equals
\[ 2\frac{(1-t^{2^{r+1}})^2}{(1-t^2)^2}; \]
\item The vector space $\coh^j(\GL,S^{2(r)})$ is abstractly isomorphic to the
vector space of coinvariants under the $\Sym_2$-action of
$\coh^j(\GL,\otimes^{2{(r)}})$; the Poincar\'e series of
$\coh_\Po^*(\GL,S^{2{(r)}})$ is equal to:
\[ \frac{(1-t^{2^{r+1}})^2}{(1-t^2)^2}+\frac{1-t^{2^{r+2}}}{1-t^4}; \]
%\item $\coh^j(\GL,S^{2(r)})~ \cong ~ S^2(E_r)^j $, any $j \geq 0$,
%where $S^2(E_r)^j$ is the degree $j$ component of the graded group $S^2(E_r)$.
\item The Poincar\'e series of $\coh_\Po^*(\GL,\Lambda^{2(r)})$ is equal to
\newline
\[ \frac{(1-t^{2^{r+1}})^2}{(1-t^2)^2}+t\frac{1-t^{2^{r+2}}}{1-t^4}; \]
\item The Poincar\'e series of $\coh_\Po^*(\GL,\Gamma^{2(r)})$ is equal to
\newline
\[ \frac{(1-t^{2^{r+1}})^2}{(1-t^2)^2}+(1+t+t^2)\frac{1-t^{2^{r+2}}}{1-t^4}. \]
\end{itemize}
\end{thm}

\begin{rem}
It is interesting that our computations show that
$\coh_\Po^*(\GL,\Lambda^{2(r)})$ can be non-zero in odd degrees (for
$p=2$), so that this is not given as $\Lambda^2(E_r)$ as one might
hope.
\end{rem}

%%%%%%%%%%%%%%%%5.1%%%%%%%%%%%%%%%%%

\subsection{Complexes and hypercohomology spectral sequences}
\label{hyper} \vskip .1in Adapting the techniques of \cite{FLS}, we
utilize the De Rham and Koszul complexes of strict polynomial
functors of degree 2, as well as the ``Symmetric complex" which is
special to characteristic 2:
\[ \xymatrix{
\Om ,d : &S^2\ar[r]^{d^1}&\otimes^2\ar[r]^{d^2}&\Lambda^2\ar[r]&0,\\
\K ,\kappa : &\Lambda^2 \ar[r]^{\kappa^1}&\otimes^2\ar[r]^{\kappa^2}& S^2\ar[r]&0,\\
S,\partial: & S^2\ar[r]^{d^1} &\otimes^2\ar[r]^{\kappa^2} &
S^2\ar[r]&0 . } \] The Koszul complex is acyclic;  the De Rham
complex has cohomology $I^{(1)}$ in degree $0$ and $1$; the
Symmetric complex has cohomology $I^{(1)}$ in degree 0.

The maps in these complexes belong to a commutative diagram:
\begin{equation}
\label{diamond}
\xymatrix{
&&S^2\ar[rrd]^{d^1}\ar[dd]&&\\
\otimes^2\ar[rru]^{\kappa^2}\ar[rrd]^{d^2}&&&&\otimes^2\\
&&\Lambda^2\ar[rru]^{\kappa^1}&&.
}
\end{equation}
Here, if $\tau$ denotes the permutation of the two factors in
$\otimes^2$, the composite from left to right is the norm map
$1+\tau$.

We shall employ four of the  hypercohomology spectral sequences
obtained by applying $\Hom_{\B}(\Gamma^{2^{r+1}}\g,-)$ to injective
resolutions of the De Rham, Koszul, and Symmetric complexes.
Namely,
\begin{itemize}
\item ``1st De Rham" ${}_{\Om}\I:$\ \ $E_1$ page $ {}_{\Om}\I_1$ has the form
\[  \coh_\Po^*(\GL ,S^{2(r)})\stackrel{d_*^1}\to
\coh_\Po^*(\GL ,\otimes^{2 (r)})\stackrel{d_*^2}\to\coh_\Po^*(\GL
,\Lambda^{2(r)}). \]

\item ``2nd De Rham" ${}_{\Om}\II$:\ \  $E_2$ page $ {}_{\Om}\II_2$ has the form
\[ \coh_\Po^*(\GL ,I^{(r)}) \stackrel{d_2^{*,1}}\to \coh_\Po^{*+2}(\GL ,I^{(r)}). \]

\item ``1st Koszul" ${}_{\K}\I$: $E_1$ page ${}_{\K}\I_1$ has the form
\[  \coh_\Po^*(\GL ,\Lambda^{2(r)})\stackrel{\kappa_*^1}\to
\coh_\Po^*(\GL ,\otimes^{2 (r)})\stackrel{\kappa_*^2}\to\coh_\Po^*(\GL ,S^{2(r)}).
 \]

\item ``1st Symmetric" ${}_S\I:$\ \ $E_1$ page $ {}_S\I_1$ has the form
\[  \coh_\Po^*(\GL ,S^{2(r)})\stackrel{d_*^1}\to
\coh_\Po^*(\GL ,\otimes^{2 (r)})\stackrel{\kappa_*^2}\to\coh_\Po^*(\GL ,S^{2(r)}).
 \]

\end{itemize}

%%%%%%%%%%%%%%%%%%%5.2%%%%%%%%%%%%%%%%%%%
\subsection{Computation of $\coh_\Po^*(\GL ,\otimes^{2 (r)})$}
\label{2step}

The special case $n = 2$ of Theorem \ref{computeI} establishes
the first assertion of Theorem \ref{deg2}

More explicitly, $\coh_\Po^*(\GL ,\otimes^{2(r)})$
is isomorphic, as a graded representation of the symmetric group $\Sym_2$, to
the $\Sym_2$-module $E_r^{\otimes 2}\otimes
k[\Sym_2]$, where the permutation $\tau$ exchanges the two tensor factor
and $E_r = \Ext_{\Po}(I^{(r)},I^{(r)})$ as in Theorem \ref{fs} . Note
that because the symmetric group $\Sym_2$ is commutative, the action by
conjugation on $k[\Sym_2]$ is trivial. So, for $0 \leq n <2^{r+2}-2$,
$\coh^n(\GL ,\otimes^{2(r)})$
is zero in odd degrees, it is free over $\Sym_2$ if $n \equiv 2~ (mod~4)$, modulo $4$, and
for $n \equiv 0 ~ (mod~4)$, it  is a sum of a free $\Sym_2$-module with two trivial
$\Sym_2$-modules generated by classes of the form $x\otimes x\otimes 1$ and
$x\otimes x\otimes \tau$.

%%%%%%%%%%%%%%%%%%%5.2%%%%%%%%%%%%%%%%%%%
\subsection{Computation of $\coh_\Po^*(\GL ,S^{2 (r)})$}
\label{Antoine} We reproduce Antoine Touz\'e's method in \cite{T}.
His argument can be summarized in degree $2$ as follows. There is a
short exact sequence due to Akin, Buschbaum and Weyman \cite[Theorem
III.1.4]{ABW}:
\[
0 \to \hc(\Lambda^2,\Lambda^2)\to S^2(gl) \to \hc(\Gamma^2,S^2)\to 0
\]
which can be readily checked by direct inspection.
After taking $r$ twists and  applying
$\Ext_\Po^*(\Gamma^{2^{r+1}},-)$, this reduces the computation of
$\coh_\Po^*(\GL ,S^{2 (r)})$ to the separable functor case. Here,
the relevant cohomology $\coh_\Po^*(\GL ,\hc(\Lambda^{2
(r)},\Lambda^{2 (r)}))\cong\Ext_\Po^*(\Lambda^{2 (r)},\Lambda^{2
(r)})$ and $\coh_\Po^*(\GL ,\hc(\Gamma^{2 (r)},S^{2
(r)}))\cong\Ext_\Po^*(\Gamma^{2 (r)},S^{2 (r)})$ are known
\cite[Theorem 5.8]{FFSS}. They are concentrated in even degrees, so
that the cohomology long exact sequence gives a vector space
isomorphism:
\[
\coh_\Po^*(\GL ,S^{2 (r)})\cong\Ext_\Po^*(\Gamma^{2 (r)},S^{2
(r)})\oplus\Ext_\Po^*(\Lambda^{2 (r)},\Lambda^{2 (r)}).
\]
The second assertion of Theorem \ref{deg2} follows.

%%%%%%%%%%%%%%%%%%%5.4%%%%%%%%%%%%%%%%%%%%%%

\subsection{Partial determination of  ranks of $d_*^1, d_*^2, \kappa_*^1, \kappa_*^2$}
\vskip .1in We proceed to consider the differentials of the first
page of the first hypercohomology De Rham and Koszul spectral
sequences.

The differentials on the first page of the first hypercohomology
spectral sequences are  the maps induced by $d$ and $\kappa$. Upon
applying $\Ext_\Po^*(\Gamma^{2^{r+1}},-)$ to (\ref{diamond}), we
obtain the commutative diagram:
\begin{equation}
\label{diamond2}
\xymatrix{
& \coh_\Po^*(\GL ,S^{2(r)})\ar[rd]^{d_*^1}\ar[dd]&\\
\coh_\Po^*(\GL ,\otimes^{2 (r)})
    \ar[ru]^{\kappa_*^2}\ar[rd]^{d_*^2}
    &&\coh_\Po^*(\GL ,\otimes^{2 (r)})\\
&\coh_\Po^*(\GL ,\Lambda^{2(r)})\ar[ru]^{\kappa_*^1}&.
}
\end{equation}
The composite from left to right is induced by the norm map
$1+\tau$. From the determination of $\coh_\Po^*(\GL
,\otimes^{2(r)})$ made explicit in \S \ref{2step}, it follows that
\begin{equation}
\label{ineq1} \rk((1+\tau)_*) = \ker((1+\tau)_*) =
\dim_n(E_r^{\otimes 2}), \quad n \equiv 2 ~(mod~4)
\end{equation}
\begin{equation}
\rk((1+\tau)_*) =  \dim_n(E_r^{\otimes 2})-1, \quad \ker((1+\tau)_*)
=  \dim_n(E_r^{\otimes 2})+1, \quad n \equiv 0 ~(mod~4)
\end{equation}
for $0 \leq n < 2^{r+2}-2$, where $\rk(-)$ denotes the dimension of
the image, $\ker(-)$ denotes the dimension of the kernel, and
$\dim_n(-)$ denotes the dimension of the homogeneous part of degree
$n$.

The diagram (\ref{diamond2}) implies the following inequalities.
\begin{equation}\label{rank inequalities}
\rk((1+\tau)_*)=\rk (d^1_*\kappa^2_*)\leq \rk (d^1_*)\leq\rk
(\kappa^1_*) \leq \ker( \kappa^2_*) \leq
\ker(d^2_*)\leq\ker((1+\tau)_* ).
\end{equation}

Thus, (\ref{ineq1}) implies that
\begin{equation}\label{2 mod 4:eq}
\rk (d^1_*)~=~ \rk (\kappa^1_*) ~=~ \ker (\kappa^2_*)~=~ \rk
(\kappa^2_*)  ~=~ \ker (d^2_*)~=~ \rk (d^2_*) ~ = ~ \dim
(E_r^{\otimes 2})
\end{equation}
in degrees congruent to 2 modulo 4 and that
\[ \dim(E_r^{\otimes 2})+1 \geq \rk(\kappa^2_*) ~ \geq ~ \dim(E_r^{\otimes 2})-1 \]
in degrees congruent to 0 modulo 4.

\begin{prop}
\label{1van} The columns ${}_K I_2^{1,*}, ~ {}_S I_2^{1,*}$ of the
$E_2$-pages of the 1st Koszul and Symmetric spectral sequences have
all entries equal to 0. The column ${}_\Om I_2^{1,*}$ of the
$E_2$-pages of the 1st De Rham spectral sequence is at most
one-dimensional in each degree $*$.
\end{prop}

\begin{proof}
Recall that the columns ${}_\Om I_1^{1,*}, ~ {}_K I_1^{1,*}, ~ {}_S
I_1^{1,*}$ are each given by $E_r^{\otimes 2} \otimes k\Sym_2$ which
vanishes in odd degree. Because the Koszul spectral sequence
converges to 0, and because the Symmetric spectral sequence
converges to $E_{r+1}$ which vanishes in even degrees, we conclude
that ${}_K I_2^{1,*} ~ = ~ {}_S I_2^{1,*} ~ = ~ 0.$ Since the second
page of the second De Rham spectral sequences is one-dimensional in
each degree, we conclude that ${}_\Om I_2^{1,n}$ itself is at most
one-dimensional.
\end{proof}

\begin{prop}
\label{2eq}
We have the following equalities among the ranks of the $d_1$ differentials in
the first hypercohomology spectral sequences:
\begin{equation}
\label{eq1} \rk (\kappa^1_*) ~ = ~\rk (d^1_*) , \quad \rk
(\kappa^2_*) ~ = ~ \rk (d^2_*)+\epsilon_*,
\end{equation}
where $\epsilon_*$ is $0$ or $1$.
\end{prop}

\begin{proof}
Proposition \ref{1van} implies that each of the following sums
equals $2\cdot \dim(E_r^{\otimes 2})$:
\begin{equation}
\rk (\kappa^1_*)+\rk (\kappa^2_*)~=~\rk (d^1_*)+\rk
(d^2_*)+\epsilon_*~= ~\rk(d^1_*) + \rk(\kappa^2_*).
\end{equation}
This immediately implies the asserted equalities.
\end{proof}

%%%%%%%%%%%%%%%%%5.5%%%%%%%%%%%%%%%%%%%%%%%
\subsection{Determination of the differential of $ {}_{\Om}\II_2$}
    Recall that the 2nd De Rham hypercohomology spectral sequence has $E_2$-page
with only two non-zero rows:  ${}_{\Om}\II_2^{*,1},
{}_{\Om}\II_2^{*,0}$. Both these rows are given by $\coh^*(\GL,
I^{(r)}) \cong E_r$, which is 1 dimensional in non-negative even
degrees less than $2^{r+1}$ and 0 otherwise.  Thus, the following
proposition completely determines this spectral sequence.

\begin{prop}
\label{Dabut} The differential in the second page of the 2nd De Rham
hypercohomology spectral sequence
\[ d_2^{*,1}: {}_{\Om}\II_2^{*,1}  = \coh^{*}(\GL,I^{(r)}) ~ \to ~
\coh^{*+2}(\GL,I^{(r)}) =  {}_{\Om}\II_2^{*+2,0} \] is zero.

Consequently, ${}_{\Om}\II_{\infty}^{n,i}$ is one-dimensional in
every degree $n$.
\end{prop}

\begin{proof}
By general arguments on hypercohomology spectral sequences, the
differential $d_2$ is left Yoneda product by a class $e$ in
$\Ext_{\B}^2(\gl^{(r)},\gl^{(r)})$. The class $e$ is obtained by
precomposing with $\gl$ the class $e_1 \in
\Ext_{\Po}^2(I^{(r)},I^{(r)}) = \coh_\Po^2(\GL ,I^{(r)})$, where
$e_1$ is constructed in \cite{FS} as the class multiplication by
which gives the differential $d_2$ in a corresponding
hypercohomology spectral sequence.

The first row $\coh^*(\GL,I^{(r)})=\Ext_{\B}^*(\Gamma^{p^r} gl,
gl^{(r)})$ is, in degree 0, one-dimensional generated by
$\Phi^r\circ gl$, where $\Phi^r:\Gamma^{p^r} \to I^{(r)}$ is the
iterated Frobenius map. Therefore:
\[
d_2(\Phi^r\circ gl)=(e_1\circ gl)\smile(\Phi^r\circ gl)=(e_1
\Phi^r)\circ gl .
\]
This is zero because $e_1 \Phi^r =0$ for $r>0$.

 We now establish the vanishing of the $d_2$-differential in all degrees. 
Using the identification
\[
\Ext_{\B}^*(\gl^{(r)},\gl^{(r)})\cong\Ext_{\Po}^*(I^{(r)},I^{(r)})^{\otimes
2}
\]
of (\ref{identification:eq}) 
(obtained by taking an injective resolution $\hc(P_\bu,I^\bu)$ of
$\gl^{(r)}$), we write:
\[
e_1\circ \gl=\alpha e_1\otimes 1 + \beta 1 \otimes e_1
\]
for coefficients $\alpha$ and $\beta$.
We make explicit the Yoneda product
 \[
\Ext_{\B}^*(\gl^{(r)},\gl^{(r)})\otimes  \Ext_\B^*( \Gamma^d
\gl,\gl^{(r)}) \to  \Ext_\B^*( \Gamma^d \gl,\gl^{(r)})
\]
in terms of the above identifications as follows. For a class
$a\otimes b$ in $\Ext_{\Po}^*(I^{(r)},I^{(r)})^{\otimes 2}\cong
\Ext_{\B}^*(\gl^{(r)},\gl^{(r)})$ and a class $x$ in
$\Ext_{\Po}^*(I^{(r)},I^{(r)})=E_r$, this Yoneda product is given by
\[
(a\otimes b)x=axb.
\]
Since $E_r$ is a commutative algebra, we conclude
\[
d_2(x)=(e_1\circ \gl)x=(\alpha e_1\otimes 1 +\beta 1 \otimes e_1)x=
\alpha e_1 x +\beta x e_1=(\alpha+\beta) e_1 x .\] The degree $0$
computation shows that $\alpha+\beta=0$, and therefore $d_2$ is zero
in all degrees.

\end{proof}

%%%%%%%%%%%%%%%%%%%%5.6%%%%%%%%%%%%%%%%%

\subsection{1st hypercohomology spectral sequences, $E_1$-pages}
We are now able to analyze the De Rham, Koszul, and Symmetric 1st
hypercohomology spectral sequences by considering cohomology degrees
modulo four.  As seen in subsection \ref{hyper}, these have the
following appearance.  The ranks of various maps in these diagrams
are indicated with notation above  the arrows.

\begin{equation}
\label{derham}
\xymatrix @!0 @C=3.5cm @R=.8cm {
\coh_\Po^*(\GL,S^{2(r)}) \ar[r]^{d^1_*} & \coh_\Po^*(\GL,\otimes^{2(r)}) \ar[r]^{d^2_*}
& \coh_\Po^*(\GL,\Lambda^{2(r)}) \\
0 & 0 & \lambda_{4n+3} \\
s_{4n+2}=g_n \ar[r]^{=g_n} & 2g_n \ar[r]^{=g_n} & \lambda_{4n+2} \\
0 & 0  & \lambda_{4n+1} \\
s_{4n}=f_n+1 \ar[r] & 2f_n \ar[r] & \lambda_{4n}}
\end{equation}

\begin{equation}
\label{koszul}
\xymatrix @!0 @C=3.5cm @R=.8cm {
\coh_\Po^*(\GL,\Lambda^{2(r)}) \ar[r]^{\kappa^1_*} & \coh_\Po^*(\GL,\otimes^{2(r)})
\ar[r]^{\kappa^2_*}
& \coh_\Po^*(\GL,S^{2(r)}) \\
\lambda_{4n+3} & 0 & 0 \\
\lambda_{4n+2} \ar[r]^{=g_n} & 2g_n \ar[r]^{=g_n} & s_{4n+2}=g_n \\
\lambda_{4n+1} & 0  & 0 \\
\lambda_{4n} \ar[r] & 2f_n \ar[r] & s_{4n}=f_n+1}
\end{equation}

and

\begin{equation}
\label{symmetric}
\xymatrix @!0 @C=3.5cm @R=.8cm {
\coh_\Po^*(\GL,S^{2(r)}) \ar[r]^{d^1_*} & \coh_\Po^*(\GL,\otimes^{2(r)}) \ar[r]^{\kappa^2_*}
& \coh_\Po^*(\GL,S^{2(r)}) \\
0 & 0 & 0 \\
s_{4n+2}=g_n \ar[r]^{=g_n} & 2g_n \ar[r]^{=g_n} & s_{4n+2}=g_n \\
0 & 0  & 0 \\
s_{4n}=f_n+1 \ar[r] & 2f_n \ar[r] & s_{4n}=f_n+1}
\end{equation}
where $n$ is some cohomological degree $< 2^{p^{r-1}}$, $2f_n =
\dim_k \coh^{4n}(\GL,\otimes^{2(r)}), 2g_n = \dim_k
\coh^{4n+2}(\GL,\otimes^{2(r)})$, where $s_m, \lambda_m$ are the
dimensions of the indicated cohomology groups (e.g., $s_{4n} =
\coh^{4n}(\GL,S^{2(r)})$), and where the equalities above the arrows
indicate the dimensions of the ranks of the maps. The computation of
$s_*$ is the second part of Theorem \ref{deg2} (see Section
\ref{Antoine}).

\begin{prop}\label{0 mod 4}
In degrees congruent to $0$ modulo $4$:
\[
\rk (d^1_*)~=~ \rk (\kappa^1_*) ~=~  \rk (\kappa^2_*)  ~=~ \rk
(d^2_*)+1 ~ = ~ \dim (E_r^{\otimes 2})
\]
\end{prop}
\begin{proof}
Looking at (\ref{symmetric}), whose abutment is $E_r$, we conclude
that in degree $4n$
\[
\rk (d^1_*)~=~  \rk (\kappa^2_*) ~ = ~ \dim (E_r^{\otimes 2})=f_n.
\]
This allows to compute the 1st Koszul spectral sequence
(\ref{koszul}). It readily implies:
\[
\begin{aligned}
\rk (\kappa^1_*) ~ &= ~ \dim (E_r^{\otimes 2})~=~f_n\\
\lambda_{4n}~&=~ f_n\\
\lambda_{4n+1}~&=~ 1\\
\lambda_{4n+2}~&=~ g_n\\
\lambda_{4n+3}~&=~ 0 ,
\end{aligned}
\]
proving the third part of Theorem \ref{deg2} as well. Finally,
looking at the 1st De Rham spectral sequence and its abutment in
Proposition \ref{Dabut}, one concludes that: $\rk (d^2_*)~=~f_n-1$
in degrees congruent to $0$ modulo $4$.
\end{proof}

%%%%%%%%%%%%%%5.8%%%%%%%%%%%%%

\subsection{Completion of proof of Theorem \ref{deg2}}
The computation in the proof of Proposition \ref{0 mod 4} explicitly
asserts that the dimension of $\coh^*_\Po(GL,\Lambda^{2(r)})$ is
half the dimension of $\coh^*_\Po(GL,\otimes^{2(r)})$, plus one
extra dimension in every degree congruent to $1$ modulo $4$. This is
the third part of Theorem \ref{deg2}.

To obtain the Poincar\'e series of $\coh_\Po^*(\GL,\Gamma^{2(r)})$,
we use the dual of the (twisted, exact) Koszul complex:
\[
0 \to~ \Gamma^{2(r)} ~ \to ~ \otimes^{2(r)} ~ \stackrel{d^2}\to ~
\Lambda^{2(r)} ~\to 0 .
 \]
By (\ref{2 mod 4:eq}) and Proposition \ref{0 mod 4}, we know in
every degree the map induced by $d^2$ in cohomology. The long exact
sequence in cohomology then implies that
$\coh_\Po^*(\GL,\Gamma^{2(r)})$ has half the dimension of
$\coh^*_\Po(GL,\otimes^{2(r)})$, plus one extra dimension in every
degree congruent to $0$, $1$ or $2$ modulo $4$, through degree $\leq
2^{r+2}-2$.

%%%%%%%%%%%%%%%%%%%%%%%%%%%%%%%
%Section 6
%%%%%%%%%%%%%%%%%%%%%%%%%%%%%%%%

\section{Computations of cohomology for $A^d(\gl^{(r)}) \otimes \gl^{(r)} $}

We recall that a (strict polynomial) functor $A $ in $ \Po_d$ is
said to be of exponential type (cf. \cite{FFSS}) if it belongs to a
sequence of functors $A^0,\ldots,A^d = A, A^{d+1} \ldots $ such that
$A^n(V\oplus W) = \bigoplus_{i+j = n} A^i(V)\otimes A^j(W)$.
Necessarily, $A^0 = 0$; we shall assume that $A^1 ~ = ~ I$ (i.e, the
identity functor in $\Po_1$).

\begin{thm}
\label{type}
Let $ A^d \in \Po_d$ be of exponential type for some $d > 1$.  Then
\begin{equation}
\label{exp} \coh_\Po^*(\GL,A^d(\gl^{(r)}) \otimes \gl^{(r)}) ~ \cong
~ \left(\coh_\Po^*(\GL,A^d(\gl^{(r)})) \oplus
\coh_\Po^*(\GL,A^{d-1}(\gl^{(r)}) \otimes \gl^{(r)})\right) \otimes
E_r.
\end{equation}

\end{thm}

\begin{proof}
We repeatedly use the adjunction isomorphisms of Proposition \ref{adjunct}.

\[
\Ext _{\B}^*(\Gamma^{(d+1)p^r}\gl,A^d(\gl^{(r)}) \otimes \gl^{(r)})
:=  \]
\[
\Ext _{\B}^*(\Gamma^{(d+1)p^r}\hc(-_1,-_2), A^d(\hc(-_1,-_2)^{(r)}
\otimes \hc(-_1,-_2)^{(r)}) =  \]
\[
\Ext _{\Po^{op} \times \Po \times
\Po^{op}}^*(\Gamma^{(d+1)p^r}\hc(-_1,-_2), A^d\hc(-_1,-_2)^{(r)}
\otimes \hc(-_3,-_2)^{(r)} \circ D_{1,3}) \] which by adjunction
equals
\[
\Ext _{\Po^{op} \times \Po \times
\Po^{op}}^*(\Gamma^{(d+1)p^r}\hc(-_1+-_3,-_2), A^d\hc(-_1,-_2)^{(r)}
\otimes \hc(-_3,-_2)^{(r)}.) \]

    In order to compactify the notation, we replace
$\Ext _{\Po^{op} \times \cdots \times \Po}^*(-,-)$ by $[-,-]$ and we
replace $\hc$ by $\h$. Expanding out the exponential functor
$\Gamma^{(d+1)p^r}(V\oplus W)$ and dropping terms with degrees in
some variable that do not match in contravariant/covariant entries, we
verify that the previous line equals in our new notation
\[
[\Gamma^{dp^r}\h(-_1,-_2) \otimes \Gamma^{p^r}\h(-_3,-_2),
A^d\h(-_1,-_2)^{(r)} \otimes \h(-_3,-_2)^{(r)}]. \] We apply
adjunction again to get that this equals
\[
[\Gamma^{dp^r}\h(-_1,-_2) \otimes \Gamma^{p^r}\h(-_3,-_4),
A^d\h(-_1,-_2+-_4)^{(r)} \otimes \h(-_3,-_2+-_4)^{(r)}]. \] We now
expand $A^d\h(-_1,-_2+-_4)^{(r)} \otimes \h(-_3,-_2+-_4)^{(r)}$
using the exponential property for $A^*$ and the additivity of
$\hc$.  By once again dropping the terms with degrees in some
variable that do not match in contravariant/covariant entries, we obtain
\[
[\Gamma^{dp^r}\h(-_1,-_2) \otimes \Gamma^{p^r}\h(-_3,-_4),
A^d\h(-_1,-_2)^{(r)} \otimes \h(-_3,-_4)^{(r)}] ~ \oplus \]
\[
[\Gamma^{dp^r}\h(-_1,-_2) \otimes \Gamma^{p^r}\h(-_3,-_4),
A^{d-1}\h(-_1,-_2)^{(r)} \otimes \h(-_1,-_4) \otimes \h(-_3,-_2)].
\] After rearranging variables on the right hand side, we apply the
K\'unneth theorem to obtain
\[
[\Gamma^{dp^r}\h(-_1,-_2),A^d\h(-_1,-_2)^{(r)}]  \otimes
[\Gamma^{p^r}\h(-_3,-_4),\h(-_3,-_4)^{(r)} ] ~ \oplus \]
\[
[\Gamma^{dp^r}\h(-_1,-_2),A^{d-1}\h(-_1,-_2)^{(r)} \otimes
\h(-_1,-_2)^{(r)}] \otimes
[\Gamma^{p^r}\h(-_3,-_4),\h(-_3,-_4)^{(r)} ]. \] This is the
computation asserted in (\ref{exp}), since $E_r =
[\Gamma^{p^r}\hc(-,-),\hc(-,-)^{(r)}].$
\end{proof}

Combining Theorems \ref{deg2} and \ref{type}, we obtain the following explicit
computations for $k$ a field of characteristic 2.

\begin{cor}
If $k$ is a field of characteristic 2 and $r \geq 0$ a non-negative integer, then
\[ \begin{aligned}
\coh_\Po^*(\GL,S^{2(r)} \otimes I^{(r)}) ~ \cong ~
(\coh_\Po^*(\GL,S^{2(r)}) \otimes E_r) \oplus (E_r^{\otimes 3}
\oplus E_r^{\otimes 3}) \\
\coh_\Po^*(\GL,\Lambda^{2(r)} \otimes I^{(r)}) ~ \cong ~
(\coh_\Po^*(\GL,\Lambda^{2(r)}) \otimes E_r)
\oplus (E_r^{\otimes 3} \oplus E_r^{\otimes 3})\\
\coh_\Po^*(\GL,\Gamma^{2(r)} \otimes I^{(r)}) ~ \cong ~
(\coh_\Po^*(\GL,\Gamma^{2(r)}) \otimes E_r) \oplus (E_r^{\otimes 3}
\oplus E_r^{\otimes 3}).
\end{aligned} \]
%If $A^2$ denotes $S^2, \Lambda^2$ or $\Gamma^2$ and if the maps
%\[ \coh_\Po^*(\GL,\otimes^{3(r)}) \to \coh_\Po^*(\GL,A^{2(r)}\otimes I^{(r)})
% \to \coh_\Po^*(\GL,\otimes^{3(r)}) \]
%are induced by the natural maps $\otimes^2 \to A^2 \to \otimes^2$, then these maps
%are given as the sum of the maps identified in the previous section,
%$\coh_\Po^*(\GL,\otimes^{2(r)}) \to \coh_\Po^*(\GL,A^{2(r)}) \to
%\coh_\Po^*(\GL,\otimes^{2(r)})$, tensored with $E_r$, plus the inclusion/projection
% of two of remaining four summands of
%$\coh_\Po^*(\GL,\otimes^{3(r)}) \simeq E_r^{\otimes 3} \otimes k\Sym_3$.
\end{cor}

\vskip .5in

%%%%%%%%%%%%%%%%%%%%%%%%%%%%%%%%
% Section 7
%%%%%%%%%%%%%%%%%%%%%%%%%%%%%%%%%

\section{Comparison with cohomology of  bifunctors}

    Throughout this section, our base field $k$ will be assumed finite
of order $q = p^e$ for some prime $p$.
We now consider the category $\F$ of all functors from the category
$\V$ of finite dimensional vector spaces over our base field $k$ to the
category of all $k$ vector spaces, and the category $\F^{op} \times \F$
of bifunctors.  We wish to study the  induced map on cohomology of the
forgetful functors
\[ \Po ~ \to ~ \F, \quad  \quad \Po^{op} \times \Po~ \to ~ \F^{op} \times \F. \]
A significant aspect of the forgetful functor $\Po \to \F$ is that
the Frobenius twist $I^{(e)} \in \Po$ maps to the identity functor
in $\F$ and thus each $I^{(r)} $ in $ \Po$ becomes an invertible
element (with respect to composition) in $\F$.

Recall that $F \in \F$ is said to be {\it finite} if it takes values
in $\V$ (i.e., $F(V)$ is finite dimensional for every finite
dimensional vector space $V$) and if it has finite Eilenberg-MacLane
degree (i.e., for some $n$, the $n$-th difference functor
$\Delta^nF$ is 0). Any functor in $\F$ which is a strict polynomial
functor (i.e., in the image of the forgetful functor $\Po \to \F$)
is finite. As shown in \cite{S}, \cite{FLS}, any finite functor has
a resolution by projective objects in $\F$ which are finite direct
sums of functors of the form $k[\Hom_k(W,-)]$ for some $W $ in $ \V$
(defined by sending $V \in \V$ to the $k$-vector space on the
underlying set of $\Hom_k(W,V)$).

 For any $F $ in $ \F$ and any $i, 1 \leq i \leq q-1$,
 we define $F^i \subset F$ to be the
subfunctor whose value on $V \in \V$ consists of those elements $x \in F(V)$
such that $F(\mu): F(V) \to F(V)$ maps $x$ to $ \mu^i x$ for any $\mu \in k^*$.
If $\phi: F \to G$ is a map in $\F$, then  the map $\phi$
restricts to a map $F^i \to G^i$ for each $i$.  If $F$ is
a finite functor, then $F \cong \oplus_i F^i$.  If $P \in \Po_d$ is a strict
polynomial functor of degree $d$, then its image under the forgetful
functor is a finite functor of weight $d$ (where $d$ is taken modulo $q-1$).

We obtain injective objects in $\F$ by dualization of projectives:
the injective $(k[\Hom_k((-)^\#,W)]) \in \F$ sends $V \in \V$ to the
$k$-vector space on the underlying set of $\Hom_k(V^\#,W)$, where
$V^\#$ is the linear dual of $V$.  Thus, we have basic projective
and injective objects in $\F^{op} \times \F$,
\[
P_{W_1,W_2} = k[\Hom_k[(-),W_1)] \otimes k[\Hom_k(W_2,(-)] \]
\[  I_{W_1,W_2} =
k[\Hom_k(W_1,(-)^\#)] \otimes k[\Hom_k((-)^\#,W_2)].
 \]

The following somewhat ad hoc definition will be adequate for our purposes.

\begin{defn}
A bifunctor $F $ in $ \F$ is said to be finite if it satisfies the following
three conditions:
\begin{enumerate}
\item $F(V,W)$ is
finite dimensional for all $V,W $ in $ \V$.
\item $F$ admits a
resolution by projectives which are finite direct sums of basic
projective bifunctors (i.e., of the form $P_{W_1,W_2}$ with $W_1,W_2
$ in $ \V$); and
\item $\F$ admits a resolution injectives which are
finite direct sums of basic injective bifunctors (i.e., of the form
$I_{W_1,W_2}$ with $W_1,W_2 $ in $ \V$).
\end{enumerate}

For example, any bifunctor of the form $\hc(A_1,A_2)$
is finite whenever $A_1,A_2 \in \F$ are finite functors.
\end{defn}

\begin{prop}
The forgetful functor $\Phi: \B \to \F^{op} \times \F$ sends a strict
polynomial bifunctor of bounded degree to a finite bifunctor.  In other
words, any strict polynomial bifunctor $T$ admits a resolution by projectives
which are finite direct sums of basic projective bifunctors $P_{W_1,W_2}$ and a
resolution by injectives which are finite direct sums of basic injective bifunctors
$I_{W_1,W_2}$.
\end{prop}

\begin{proof}
We establish the existence of such a projective resolution for any strict
polynomial bifunctor $T$; the argument for the existence of a corresponding
injective resolution is similar.

We first assume that $T = \hc(A,B)$ is of separable type, where $A, B$ are
strict polynomial functors.  Since every strict polynomial functor $P$ is a finite
functor (as observed in \cite{FFSS}), a resolution of $A$ by finite direct sums of
basic injective functors and a resolution of $B$ by finite direct sums of basic
projective functors determine (upon external tensor product) a resolution of $T$
by finite direct sums of basic projective bifunctors.

A general strict polynomial bifunctor $T$ admits a resolution by strict polynomial
functors of separable type by Proposition \ref{enuf}.  Thus, taking a resolution of
each term in such a resolution, we obtain a bicomplex consisting of finite
direct sums of basic projective bifunctors, and its total complex gives us the
required resolution of $T$.
\end{proof}

    A special role is played by the bifunctor
\[ k[\gl] =: ~ k[\Hom_k(-,-)] ~ \in \F^{op} \times \F \]
as we first see in the following analogue of (\ref{sep}).
This is proved in a manner exactly parallel to the proof
of Proposition \ref{spec}.

\begin{prop}
\label{spec2} For two finite functors $A_1$ and $A_2$, we consider the
bifunctor $\hc (A_1,A_2)$ of separable type.  For any bifunctor $F$,
there is a convergent spectral sequence of the form
\begin{equation}
\label{spec3} E_2^{s,t} ~ = ~
\Ext_{\F}^s(A_1(-_1),\Ext_{\F}^t(F(-_1,-_2),A_2(-_2))) ~ \Rightarrow
~ \Ext_{\F^{op} \times \F}^{s+t}(F,\hc (A_1,A_2)),
\end{equation}
natural in $A_1$, $A_2$ and $F$.

If $F = k[\gl] $, then this spectral sequence collapses to give the
natural isomorphism
\begin{equation}
\label{collapse} \Ext_{\F^{op} \times \F}^*(k[\gl],\hc (A_1,A_2))
\cong  \Ext_{\F}^*(A_1,A_2).
\end{equation}

If $F = \hc(B_1,B_2)$ with $B_1, B_2$ also finite functors, then this spectral
sequence collapses to give the natural isomorphism
\begin{equation}
\label{collapse2}
\Ext_{\F^{op} \times \F}^*(\hc(B_1,B_2),\hc (A_1,A_2)) \cong
\Ext_{\F}^*(B_1,B_2) \otimes \Ext_{\F}^*(A_1,A_2).
\end{equation}

\end{prop}

\begin{proof}
We first verify by inspection that
\[ \Hom_{\F^{op} \times \F}(P_{W_1,W_2},F) ~ \cong ~ F(W_1,W_2) ~ \cong ~
\Hom_{\F^op \times \F}(F,I_{W_1,W_2}). \] As in Lemma \ref{switch},
this implies the isomorphism
\[
\Hom_\biF(F(-_1,-_2),k[\Hom_k(W_1,(-_1)^\#)] \otimes
k[\Hom_k((-_2)^\#,W_2)]) ~ \cong ~ \]
\[ \Hom_\F(\Hom_\F(F(-_1,-_2),k[\Hom_k(W_1,(-_1)^\#)]),k[\Hom_k((-_2)^\#,W_2)]). \]
Taking $F = k[\gl]$, we proceed exactly as in the proof of
Proposition \ref{spec} to establish the spectral sequence
(\ref{spec3}) and prove (\ref{collapse}).

To prove (\ref{collapse2}), we argue exactly as in the proof of Proposition \ref{Ku}.
\end{proof}

Using a theorem of A. Suslin \cite[A.1]{FFSS}, we obtain the following
interpretation of bifunctor cohomology.  We introduce the notation
\[ \coh^*(\GL(k),F) ~ =: ~  \coh^*(\GL(n,k),F(k^n,k^n)), \quad n >>0 \]
for a finite bifunctor $F$, where $\coh^*(\GL(n,k),F(k^n,k^n))$
denotes the usual group cohomology of the finite group $\GL(n,k)$
with coefficients in the $\GL(n,k)$-module $F(k^n,k^n)$.

\begin{thm}
\label{relate}
Let $F \in \F^{op} \times \F$ be a finite bifunctor.  Then the natural map
\[ \Ext_{\F^{op} \times \F}^*(k[\gl],F) ~ \to ~  \coh^*(\GL(k),F) \]
is an isomorphism.

Moreover, for a given cohomological degree $s$,
\[ \Ext_{\F^{op} \times \F}^s(k[\gl],F) ~ \cong ~  \coh^s(\GL(n,k),F(k^n,k^n)),
\quad n >> s. \]

\end{thm}

\begin{proof}
 For $n >> s$ and $A_1, A_2$ finite functors, a theorem of W. Dwyer
 \cite{Dw} asserts that
 the natural map
\[ \Ext_{\GL(n,k)}^s(A_1(k^n),A_2(k^n)) ~ \to ~ \Ext_{\GL(n+1,k)}^s(A_1(k^{n+1}),A_2(k^{n+1})) \]
is an isomorphism; we denote the stable value by $
\Ext_{\GL(k)}^s(A_1,A_2)$. Suslin's theorem \cite[A.1]{FFSS} asserts
that the natural map
\[ \Ext_\F^*(A_1,A_2) ~ \to ~ \Ext_{\GL(k)}^*(A_1,A_2) \]
is an isomorphism.  (Although Suslin's theorem requires that $k$ be
finite as we also require in this section, this isomorphism has been
generalized by A. Scorichenko to an arbitrary field; we refer the
interested reader to \cite{Sch} and \cite{F-P}.) This isomorphism,
Dwyer's stability \cite{Dw}, and Proposition \ref{spec2} imply the
two assertions of the theorem in the special case in which $F$ is of
the form $\hc(A_1,A_2)$ with $A_1, A_2$ finite functors.

For a general finite bifunctor $F$, we first consider cohomological degree 0.
Let $0 \to F \to I^0 \to I^1$ be exact in $\biF$ with both $I^0, I^1$ injective
objects in $\biF$ which are directs sums of (finite) bifunctors of the form
$I_{W_1,W_2}$.  Then we conclude the natural isomorphism
\[ \Hom_{\F^{op} \times \F}^0(k[\gl],F) ~ \to ~  \coh^0(\GL(k),F) \]
by the left exactness of both $\Hom_\biF$ and $\coh^0$ together with
the special case verified above; the second assertion for
cohomological degree 0 also follows from left exactness.

To complete the proof, we argue by induction on the cohomological
degree and use a devissage argument.  Namely, choose a short
exact sequence
$0 \to F \to I \to G \to 0 $ in $\biF$ with $I$ a direct sum of (finite) bifunctors
of the form $I_{W_1,W_2}$.  Then the two assertions in cohomological
degree $s-1$ for $G$ imply the corresponding assertions in cohomological
degree $s$ for $F$.
\end{proof}

In parallel with the notation $\coh^*_\Po(\GL,A$) introduced in
(\ref{note}), we employ the notation
\[ \coh^*_\F(\GL,A) =: ~ \Ext_{\F^{op}\times \F}^*(k[\gl],A \circ \gl)
~ = ~ \coh^*(\GL(k),A(\gl)) \] for a bifunctor of the form $A \circ
\gl$ with $A$ a functor of finite type.

    To convey the similarities and differences of computations of
cohomology in $\Po^{op} \times \Po$ and $\F^{op} \times \F$, we state
and prove the analogue of Theorem \ref{computeI}.

\begin{prop}
\label{computeII} Let $E_\infty =: \varinjlim_r E_r$, a divided
power algebra on the infinite sequence of generators $e_1,
\ldots, e_r, \ldots $ with the degree of $e_i$ equal to $2p^{i-1}$. For
any $n$, we have a $\Sym_n$-equivariant isomorphism
\[ \coh^*_\F(\GL,\otimes^n) ~ \cong ~ \varinjlim_r \coh^*(\GL,\otimes^{n(r)})
     ~ \cong ~ E_\infty^{\otimes n}\otimes k[\Sym_n]. \]
\end{prop}

\begin{proof}
The adjunction isomorphisms of Proposition \ref{adjunct} are equally
valid with $\Po$ replaced by $\F$ (i.e., with strict polynomial
multi-functors replaced by multi-functors).

To sort out
\begin{equation}
\label{sort}
\Ext_{\F\times \cdots \times \F}^*(\otimes^{n} \circ \bigoplus,I
\boxtimes \cdots \boxtimes I),
\end{equation}
the degree argument can be replaced by a second adjunction
comparing functors in $n-1$ and $n$ variables.  Consider the
exact functors
\begin{equation}
\epsilon_0: \V^{\times n-1} \quad \stackrel{\longrightarrow}{\longleftarrow}
\quad  \V^{\times n}:  \iota
\end{equation}
which are left and right adjunct to each other, $\epsilon_0$
 given by $(V_1,\ldots,V_{n-1}) \mapsto (V_1,\ldots,V_{n-1},0)$
and $\iota$ given by $(V_1,\ldots,V_n) \mapsto
(V_1,\ldots,V_{n-1})$.  Thus, precomposition by $ \iota$ and $\epsilon_0$
determine left and right adjunct functors on functor categories
\begin{equation}
\label{missing}
(-\circ \iota): \F^{\times n-1} \quad \stackrel{\longrightarrow}{\longleftarrow}
\quad  \F^{\times n} : (-\circ \epsilon_0).
\end{equation}
Since $\boxtimes^n \circ \epsilon_0 = 0$, all the terms
in the expansion of (\ref{sort}) with a missing variable are 0.  This
enables us to drop the same terms which were dropped for
degree reasons in the proof of Theorem \ref{computeI}.

The K\"unneth Theorem remains valid for $\Ext_\F^*$
for finite functors.  Thus, the proof
of Theorem \ref{computeI} applies to prove this proposition thanks to the
fact proved in \cite{FS} that the natural map
\[ \varinjlim_r \Ext_\Po^*(I^{(r)},I^{(r}) ~ \to ~ \Ext_\F^*(I,I) \]
is an isomorphism.
\end{proof}

The following theorems extend the applicability of Proposition
\ref{computeII} with the important restriction that the degree of
$F$ is less than or equal to the cardinality of $k$.  (An analysis
of change of field as in \cite{FFSS} can be achieved to show that
$\coh^*(\GL(K),F_K)$ is a direct factor in
$\coh^*(\GL(k),F)\otimes_k K$ for an extension $K/k$ of finite
fields.)

\begin{thm}%\label{P2F}
Let $k$ be a finite field with $q$ elements. If $S, ~ T$ are strict
polynomial bifunctors homogenous of bidegree $(d,e)$, $d,e\leq q$,
then there is a natural isomorphism:
\[ \varinjlim_r \Ext_\B^*(S^{(r)},T^{(r)}) ~ \stackrel{\simeq}{\to} ~
\Ext^*_{\F^{op}\times \F}(S,T). \]
\end{thm}

\begin{proof}
For $S = \hc(B_1,B_2), ~ T = \hc(A_1,A_2)$  of separable type, we employ
 the natural commutative square
\[ \xymatrix{
\varinjlim_r \Ext_\B^*(\hc(B_1,B_2)^{(r)}, \hc(A_1,A_2)^{(r)})\ar[r]\ar[d]&
\Ext_{\F^{op}\times \F}^*(\hc(B_1,B_2), \hc(A_1,A_2)) \ar[d]\\
\varinjlim_r \Ext_\Po^*(A_1^{(r)},B_2^{(r)}) \otimes \Ext_\Po^*(B_2^{(r)},A_2^{(r)})
\ar[r]&
\Ext_\F^*(A_1,B_2) \otimes \Ext_\F^*(B_2,A_2)
 } \]
whose vertical arrows are the isomorphisms of Propositions \ref{Ku}
and \ref{spec2}. For degrees $\leq q$, the lower horizontal arrow of
the above square is an isomorphism by \cite[3.10]{FFSS} so that the
upper horizontal arrow (the natural map in the statement of the
theorem) is also an isomorphism.

For general strict polynomial bifunctors $S, ~T$, we consider
resolutions $P_\bu \to S , ~T \to J^\bu$ of $S, ~T$ by separable
strict polynomial bifunctors of same bidegree. Then we identify the
map of the theorem with the map on hyperext-groups
\[ \varinjlim_r \bE_\B^*(P_{\bu (r)},J^{\bu (r)}) ~ \stackrel{\simeq}{\to} ~
\bE^*_{\F^{op}\times \F}(P_\bu,J^\bu). \] The proof is completed by
comparing spectral sequences for hyperext groups and employing the
special case of the theorem proved above for strict polynomial
bifunctors of separable type.
\end{proof}

\begin{thm}\label{compare} %\label{GL2GL(k)}
Let $k$ be a finite field with $q$ elements. If $T$ is a strict
polynomial bifunctor homogeneous of bidegree $(d,d)$, $d\leq q$,
then there is a natural isomorphism:
\[ \varinjlim_r \coh_\Po^*(\GL,T^{(r)})  \simeq  \coh_\F^*(\GL,T). \]
\end{thm}

\begin{proof}
For $T = \hc(A,B)$  of separable type of bidegree $(d,d)$, we employ
 the natural commutative square
\[ \xymatrix{
\varinjlim_r \Ext_\B^*(\Gamma^{p^rd}\gl,
\hc(A,B)^{(r)})\ar[r]\ar[d]&
\Ext_{\F^{op}\times \F}^*(k[\gl], \hc(A,B)) \ar[d]\\
\varinjlim_r \Ext_\Po^*(A^{(r)},B^{(r)}) \ar[r]& \Ext_\F^*(A,B)
 } \]
whose vertical arrows are isomorphisms by \ref{sep} and
\ref{collapse}. When $d\leq q$, the lower horizontal arrow of the
above square is an isomorphism by \cite[3.10]{FFSS} so that the
upper horizontal arrow (the natural map in the statement of the
theorem) is also an isomorphism. We then argue as in the proof of
the preceding theorem by resolving $T$ by separable strict
polynomial bifunctors.
\end{proof}

%%%%%%%%%%%%%%%%%%%%%%%%%%%%%%%%
% Section 8
%%%%%%%%%%%%%%%%%%%%%%%%%%%%%%%%%

\section{Consequences for cohomology of discrete groups}

    In determining the $K$-theory of $\Z/p^2\Z$ in low degrees,
Friedlander and L. Evens \cite{EF} used the Hochschild-Serre spectral sequence
for the extension
\[ 1 \to \gl(n,\bF_p) ~ \to ~ \GL(n,\Z/p^2\Z) ~ \to ~ \GL(n,\bF_p) \to 1 \]
which takes the form
\begin{equation}
\label{zp2}
E_2^{s,t} = \coh^s(\GL(n,\bF_2),S^t(\gl(n,\bF_2))) \\
 \Rightarrow \coh^{s+t}(\GL(n,\Z/p^2\Z),\bF_2).
 \end{equation}
for $\Z/2\Z$ and the form
\begin{equation}
\label{zp} E_2^{s,t} = \coh^s(\GL(n,\bF_p),\bigoplus_{t_1 + 2t_2 =
t}\Lambda^{t_1}(\gl(n,\bF_p))
\otimes S^{t_2}(\gl(n,\bF_p))) \\
 \Rightarrow \coh^{s+t}(\GL(n,\Z/p^2\Z),\bF_p).
 \end{equation}
for $p$ odd. We restate the calculations of \cite{EF} (which were
formulated in homology with $\gl$ replaced by the subfunctor of
trace 0 matrices).  In the computations which follow this statement,
we extend each of these calculations to all cohomological degrees as
well as weaken the restriction on the prime $p$.

\begin{prop}

\label{ef} (cf. \cite{EF})
Let $k$ be a finite field of characteristic $p \geq 5$ and $n \geq 2$.  Then
%\begin{itemize}
\[
\coh^i(\GL(n,k),\gl(n,k))=
\left\{%
\begin{array}{ll}
    k, & \hbox{ if i = 0, 2} \\
    0, & \hbox{ if i = 1, 3} \\
  \end{array}
\right.
 \]

\[
\coh^i(\GL(n,k),S^2\gl(n,k)) =
\left\{%
\begin{array}{ll}
    k^2 & \hbox{ if i = 0} \\
    0 & \hbox{ if i = 1} \\
  \end{array}
\right.
 \]

\[
\coh^i(\GL(n,k),\Lambda^2\gl(n,k)) =
\left\{%
\begin{array}{ll}
    0 & \hbox{ if i = 0,1} \\
    k^2 & \hbox{ if i = 2} \\
  \end{array}
\right.
 \]

\[ \coh^i(\GL(n,k),\Lambda^3\gl(n,k)) =
\left\{%
\begin{array}{ll}
    k & \hbox{ if i = 0} \\
    0 & \hbox{ if i = 1} \\
  \end{array}
\right.
 \]

\end{prop}

Our results in previous sections enable us to provide considerable
generalizations of these computations when stabilized (so that
$\GL(n,k)$ for $n \geq 2$ is replaced by $\GL(n,k)$ for $n>> 0$).
The isomorphism (\ref{finite2}) of the following theorem gives an
explicit computation in all degrees of the stabilizations of the
very low degree cohomology computed in Proposition \ref{ef}.

\begin{thm}
\label{compute}
Let $A$ be a strict polynomial functor of degree $d$ and $k$ a finite field of
order $q = p^e$ for some prime $p$.  Then
\begin{equation}
\label{finite1} \coh^s(\GL(k),A(\gl)) ~ \cong ~
\coh_\Po^*(\GL,A^{(r)}) \quad {\text if~} r \geq
log_p(\frac{s+1}{2}), ~ q \geq d.
\end{equation}
This isomorphism enables explicit computations as follows:
\begin{itemize}
\item If $p = d = 2$, then $\coh^*(\GL(k),A(\gl))$ is given explicitly by
Theorem \ref{deg2}.
\item If $p > d$ and if $\lambda$ is a partition of $d$, then
\begin{equation}
\label{finite2} \coh^s(\GL(k),S^\lambda(\gl)) ~ \cong ~
s_\lambda((E_r^{\otimes d} \otimes k\Sym_d)^s) \quad  {\text if~} r
\geq log_p(\frac{s+1}{2}).
\end{equation}
\end{itemize}
\end{thm}

\begin{proof}
Isomorphism (\ref{finite1}) follows by applying  Theorem
\ref{relate} which relates $\coh^*(\GL,-)$ to $\coh_\F^*(\GL,-)$,
Theorem \ref{compare} (requiring $d \leq q$) which relates
$\coh_\F^*(\GL,-)$ to  the``generic" cohomology
$\varinjlim_r\coh_\Po^*(\GL,(-)^{(r)})$, and  Corollary \ref{rtwist}
(requiring $r \geq log_p(\frac{s+1}{2})$) which enables one to avoid
the colimit with respect to Frobenius twisting.

Isomorphism (\ref{finite2}) follows immediately from isomorphism (\ref{finite1})
and Proposition \ref{partition} (which requires $p > d$).
\end{proof}

\begin{cor}
\label{Pseries} Let $k$ be a finite field of characteristic $p$.  We
have the following computations of the group cohomology of $\GL(k)$
with coefficients as indicated:
\begin{itemize}
\item The Poincar\'e series for $\coh^*(\GL(k),\gl)$ equals
\[ \frac{1}{1-t^2}. \]
\item The Poincar\'e series for $\coh^*(\GL(k),\gl^{\otimes 2})$ equals
\[ 2\left( \frac{1}{1-t^2} \right)^2. \]
\item The Poincar\'e series for $\coh^*(\GL(k),S^2(\gl))$ equals
\[ \frac{1}{(1-t^2)^2}
+ \frac{1}{1-t^4}. \]
\item  If $p=2$, then the Poincar\'e series for $\coh^*(\GL(k),\Lambda^2(\gl))$ equals
\[ \frac{1}{(1-t^2)^2}
+ \frac{t}{1-t^4}. \]
\item  If $p > 2$, then the Poincar\'e series for $\coh^*(\GL(k),\Lambda^2(\gl))$ equals
\[ \frac{1}{(1-t^2)^2}
- \frac{1}{1-t^4}. \]
\item  If $p=2$, then the Poincar\'e series for $\coh^*(\GL(k),\Gamma^2(\gl))$ equals
\[ \frac{1}{(1-t^2)^2}
+ \frac{1+t+t^2}{1-t^4}. \]
\end{itemize}
\end{cor}

\begin{proof}
The asserted Poincar\'e series for $p=2$ are obtained by taking the
limits with respect to $r$ of the Poincar\'e series given in Theorem
\ref{deg2}.

For $p > 2$, we apply Proposition \ref{partition}.  Observe that, in contrast
with $p=2$, for $p > 2$ there are split exact short exact sequence in cohomology
\[ 0 \to \coh^*(\GL,\Lambda^{2(r)}) \to \coh^*(\GL,\otimes^{2(r)}) \to \coh^*(\GL,S^{2(r)}) \to 0. \]
\end{proof}

The following proposition enables the computation of further terms in the
spectral sequence (\ref{zp}).

\begin{prop}
\label{typeF} Let $k$ be a finite field with $q = p^e$ elements and
consider $ A^d \in \Po_d$ of exponential type, with $d \leq q$.
Then
\begin{equation}
\label{exp in F} \coh_\F^*(\GL,A^d(\gl) \otimes \gl) ~ \cong ~
\left(\coh_\F^*(\GL,A^d(\gl)) \oplus \coh_\F^*(\GL,A^{d-1}(\gl)
\otimes \gl)\right) \otimes E_\infty.
\end{equation}

\end{prop}

\begin{proof}
The proof is very similar to that of Proposition \ref{type}; in particular, we repeatedly
use adjunction isomorphisms.
\[
\begin{aligned}
\Ext _{{\F^{op}\times \F}}^*(k[\gl],A^d(\gl)
        & \otimes \gl) :\\
    =\Ext _{{\F^{op}\times \F}}^*
        &(k[\hc(-_1,-_2)], A^d(\hc(-_1,-_2)\otimes\hc(-_1,-_2))\\
    = \Ext _{\F^{op} \times \F \times \F^{op}}^*
        &(k[\hc(-_1,-_2)],
A^d\hc(-_1,-_2) \otimes \hc(-_3,-_2) \circ D_{1,3})
\end{aligned}
 \]

which by adjunction equals
\[
\Ext _{\F^{op} \times \F \times}^*(k[\hc(-_1+-_3,-_2)],
A^d\hc(-_1,-_2) \otimes \hc(-_3,-_2)). \] Since $k[V\oplus W]\cong
k[V]\otimes k[W]$, the previous line equals:
\[
\Ext _{\F^{op} \times \F \times \F^{op}}^*(k[\hc(-_1,-_2)]\otimes
k[\hc(-_3,-_2)], A^d\hc(-_1,-_2) \otimes \hc(-_3,-_2)).
 \]

    In order to compactify the notation, we replace
$\Ext _{\F^{op} \times \cdots \times \F}^*(-,-)$ by $\{ -,-\}$ and
we replace $\hc$ by $\h$. We apply adjunction again to get that this
equals in this notation
\[
\{ k[\h(-_1,-_2)] \otimes k[\h(-_3,-_4)], A^d\h(-_1,-_2+-_4) \otimes
\h(-_3,-_2+-_4)\}.
 \]
We now expand $A^d\h(-_1,-_2+-_4)
\otimes \h(-_3,-_2+-_4)$ using the exponential property for
$A^*$ and the additivity of $\hc$. Typically, terms with missing
variable will not contribute.   In our case, for example,
\[ \begin{aligned}
{\{}k[\h(-_1,-_2)] \otimes  k[\h(-_3,-_4)],
        & A^d\h(-_1,-_2) \otimes\h(-_3,-_2)\}\\
    = \{ k[\h(-_1,-_2)] \otimes k[\h(-_3,0)],
        & A^d\h(-_1,-_2) \otimes \h(-_3,-_2)\} \\
    =\{ k[\h(-_1,-_2)],
        & A^d\h(-_1,-_2) \otimes\h(-_3,-_2)\}\\
    = \{ k[\h(-_1,-_2)],
        &A^d\h(-_1,-_2) \otimes\h(0,-_2)\} =0.
\end{aligned} \]

Using the same techniques, and $A^1=I$, one shows that
\[
\{ k[\h(-_1,-_2)] \otimes k[\h(-_3,-_4)], A^{d-1}\h(-_1,-_2) \otimes
A^1\h(-_1,-_4) \otimes\h(-_3,-_4) \} =0.
 \]
Finally, terms
\[
\{ k[\h(-_1,-_2)] \otimes k[\h(-_3,-_4)], A^{d-i}\h(-_1,-_2) \otimes
A^i\h(-_1,-_4) \otimes ...\}
 \]
vanish when $i\geq 2$ because, when $d\leq q$, their degree (modulo
$q-1$) in some variable do not match in contravariant/covariant
entries. We obtain
\[
\{ k[\h(-_1,-_2)] \otimes k[\h(-_3,-_4)], A^d\h(-_1,-_2) \otimes
\h(-_3,-_4)\} ~ + \]
\[
\{ k[\h(-_1,-_2)] \otimes k[\h(-_3,-_4)], A^{d-1}\h(-_1,-_2) \otimes
\h(-_1,-_4) \otimes \h(-_3,-_2) \}
 \]
and conclude by rearranging variables as in Proposition \ref{type}.
\end{proof}

\begin{rem}
Proposition \ref{typeF} is also valid if we replace $\gl$ with
$\gl^{(r)}$ for any $r\geq 0$.  The proof is the same, except that
many of the occurrences of $\hc(-,-)$ must be replaced by
$\hc(-,-)^{(r)}$.
\end{rem}

We combine Proposition \ref{typeF} with the computations of Corollary
\ref{Pseries} to obtain further computations.

\begin{cor}  Let $k$ be a finite field of characteristic $p$.
\begin{itemize}
\item The Poincar\'e series for $\coh^*(\GL(k),S^2(\gl)\otimes \gl)$ equals
\[ \left(\frac{3}{(1-t^2)^2} + \frac{1}{1-t^4} \right)\cdot \frac{1}{1-t^2}. \]
\item If $p=2$, the Poincar\'e series for $\coh^*(\GL(k),\Lambda^2(\gl)\otimes \gl)$
 equals
\[ \left(\frac{3}{(1-t^2)^2} + \frac{t}{1-t^4} \right)\cdot \frac{1}{1-t^2}. \]
\item If $p>2$, the Poincar\'e series for $\coh^*(\GL(k),\Lambda^2(\gl)\otimes \gl)$
 equals
\[ \left(\frac{3}{(1-t^2)^2} - \frac{1}{1-t^4} \right)\cdot \frac{1}{1-t^2}. \]
\end{itemize}

\end{cor}

We complete the computations begun in Proposition \ref{ef}
with the following computation.

\begin{prop}
Let $k$ be a finite field of characteristic $p>3$.  Then
\begin{itemize}
\item The Poincar\'e series for $\coh^*(\GL(k),\Lambda^3(\gl))$ equals
\[ \frac{1+5t^6}{(1-t^2)(1-t^4)(1-t^6)}. \]
\item The Poincar\'e series for $\coh^*(\GL(k),S^3(\gl))$ equals
\[ \frac{3+2t^2 + t^4}{(1-t^4)(1-t^6)}. \]
\end{itemize}
\end{prop}

\begin{proof}
We use Proposition \ref{partition}  in conjunction with Proposition
\ref{ps} to compute $\coh^*_\Po(\GL,\Lambda^{3(r)})$ and
$\coh^*_\Po(\GL,S^{3(r)})$. We apply Theorem \ref{compute} and take
the limit as $r$ goes to $\infty$ to obtain the assert Poincar\'e
series.
\end{proof}

%%%%%%%%%%%%%%%%%%%%%%%%%%%%
%%%%%%%%%%%%%%%%%%%%%%%%%%%%%

\end{document}